\documentclass[12pt,numreferences]{kluwer}

\usepackage{amssymb, amsmath, amsfonts}
\usepackage{bm}
\def\genfd{{\bm k}}
\def\id{\mathrm{id}}
\def\op{\mathrm{op}}
\long\def\nodo#1{{}}
\def\gg{\mathfrak{g}}

\def\AA{\mathcal{A}}
\def\CC{\mathcal{C}}
\def\OO{\mathcal{O}}
\def\mcS{\mathcal{S}}

\def\End{\operatorname{End}}
\def\Ker{\operatorname{ker}}
\def\Diff{\operatorname{Diff}}
\def\Vect{\operatorname{Vect}}
\def\hx{\hat{x}}
\def\hy{\hat{y}}
\def\ad{\operatorname{ad}}
\def\Ad{\operatorname{Ad}}

\def\MR#1{} 
\def\hx{\hat{x}}

\def\btr{\blacktriangleright}
\def\btl{\blacktriangleleft}
\newtheorem{defn}{Definition}
\newtheorem{prop}{Proposition}
\newtheorem{lem}{Lemma}
\newtheorem{theorem}{Theorem}

\begin{document}
\begin{opening}
\title{Lie algebra type noncommutative phase spaces are Hopf algebroids}
\author{{\sc Stjepan Meljanac} {\tt meljanac@irb.hr}}
\institute{
Theoretical Physics Division, Institute Rudjer Bo\v{s}kovi\'{c}
\\ Bijeni\v{c}ka cesta~54, P.O.Box 180, HR-10002 Zagreb, Croatia
}
\author{{\sc Zoran \v{S}koda} {\tt zoran.skoda@uhk.cz}}
\institute{Faculty of Science, University of Hradec Kr\'alov\'e\\
Rokitansk\'eho 62, Hradec Kr\'alov\'e, Czech Republic
\\
}
\author{{\sc Martina Stoji\'c} {\tt stojic@math.hr}}
\institute{Department of Mathematics, 
University of Zagreb\\
Bijeni\v{c}ka cesta~30, HR-10000 Zagreb, Croatia
}
\runningtitle{Hopf algebroid}
\runningauthor{Martina Stoji\'c, Zoran \v{S}koda and Stjepan Meljanac}
\begin{abstract} 
For a noncommutative configuration space whose coordinate algebra is the
universal enveloping algebra of a finite dimensional Lie algebra, 
it is known how to introduce an extension playing
the role of the corresponding noncommutative phase space, 
namely by adding the commuting deformed derivatives 
in a consistent and nontrivial way, therefore
obtaining certain deformed Heisenberg algebra. 
This algebra has been studied in physical contexts, 
mainly in the case of the kappa-Minkowski space-time. 
Here we equip the entire phase space algebra with a coproduct, 
so that it becomes an instance of a completed variant of 
a Hopf algebroid over a noncommutative base,
where the base is the enveloping algebra. 
\end{abstract}

\classification{AMS classification}{16S30,16S32, 16S35, 16Txx}

\vskip .03in
\end{opening}

\section{Introduction}

Recently, a number of physical models 
has been proposed~\cite{AC,HallidaySzabo,Luk91},
where the background geometry
is described by a noncommutative {\it configuration} space
of Lie algebra type. 
Descriptively, its coordinate algebra is the universal enveloping algebra 
$U(\gg)$ of a Lie algebra $\gg$ with basis $\hx_1,\ldots,\hx_n$ 
(noncommutative coordinates). So-called $\kappa$-Minkowski space
is the most explored example~\cite{Govstat,Luk91,LukNowHeis,kappaind}. 
That space has been used to build a model featuring the double special
relativity, a framework modifying special relativity, 
proposed to explain some phenomena observed in the 
high energy gamma ray bursts. 

The noncommutative {\it phase} space of the Lie algebra $\gg$
is introduced by enlarging $U(\gg)$ 
with additional associative algebra generators, 
the {\it deformed derivatives}, 
which act on $U(\gg)$ via an action $\blacktriangleright$ satisfying 
deformed Leibniz rules~\cite{scopr,heisd}.  
The subalgebra generated by the 
deformed derivatives is commutative.
In fact, this commutative algebra is a topological Hopf algebra
isomorphic to the full algebraic dual $U(\gg)^*$ 
of the enveloping algebra. 
In this article, we extend the coproduct of 
the topological Hopf algebra $U(\gg)^*$ 
of deformed derivatives to a coproduct 
$\Delta:H\to H\hat\otimes_{U(\gg)} H$ on
the whole phase space $H$ (and its completion $\hat{H}$); 
this coproduct is moreover a part of 
a (formally completed) Hopf algebroid structure on $H$ 
over the noncommutative base algebra $\AA = U(\gg)$. Descriptively,
a Hopf algebroid is an associative bialgebroid 
(Definition~\ref{def:leftbialg} in~Subsection~\ref{ss:algprel})
with an antipode map (the antipode is treated in Section~\ref{sec:antipode}).

The notion of a Hopf algebroid in this paper is slightly modified 
regarding that the tensor product $\hat\otimes_{U(\gg)}$ 
in the definition of the coproduct is understood in a completed sense;
a part of the definition still needs the tensor products without
completions. Our bialgebroid structure is similar but a bit weaker than
the bialgebroid {\it internal}~\cite{bohmInternal} 
to the tensor category 
of complete cofiltered vector spaces; 
a true internal variant is possible
in a more intricate monoidal category involving 
filtrations of cofiltrations and 
is treated along with generalizations in~\cite{StojicPhD}. 

The noncommutative phase space of Lie type is nontrivially
isomorphic to an infinite-dimensional version 
$U(\gg)\sharp U(\gg)^*\cong U(\gg)\sharp\hat{S}(\gg^*)$
of the Heisenberg double of $U(\gg)$~\cite{heisd}.
Heisenberg doubles of finite dimensional Hopf algebras
are known to carry a Hopf algebroid structure~\cite{BrzMilitaru,Lu}.
However, our starting Hopf algebra $U(\gg)$ is infinite-dimensional, 
though filtered by finite-dimensional pieces. 
While the generalities on such filtered algebras can be used to
obtain the Hopf algebroid structure~\cite{StojicPhD},
we here use the specific features of $U(\gg)$ instead, and in particular
the matrix $\OO$ introduced in the Section~\ref{sec:defphasesp}
and used to define the crucial part of the bialgebroid structure, 
the {\it target} map 
$\beta:\hx_\alpha\mapsto\sum_\beta\hx_\beta\otimes(\OO^{-1})^\beta_\alpha$.

From a dual geometric viewpoint, where $U(\gg)$ is 
viewed as the algebra of left invariant differential operators 
on a Lie group, the matrix $\mathcal{O}$ 
is interpreted as a transition matrix
between a basis of left invariant and a basis 
of right invariant vector fields. Then  
our phase space appears as the algebra of 
formal differential operators $\Diff^\omega(G,e)$ around the unit $e$ of
the Lie group $G$ integrating $\gg$. Every formal differential
operator is a finite sum of products of the form $f_s D_s$ where
$D_s$ is an invariant differential operator (belonging to $U(\gg)$) and 
$f_s$ is a formal function 
(this decomposition amounts to a Hopf algebraic
smash product in the algebraic part of the paper). 
By L.~Schwartz's theorem~\cite{CartierHopf}, the space $J^\infty(G,e)$
of formal functions at $e$ 
and $U(\gg)$ are dual by evaluating the differential operator
on a function at $e$; the duality equips $J^\infty(G,e)$ 
with a topological Hopf algebra structure with coproduct $\Delta^J$. 
Then define the coproduct of the $J^\infty(G,e)$-bialgebroid 
$\Diff^\omega(G,e)$ by the {\it scalar extension} formula 
$\Delta(f_s D_s) = (f_s\otimes 1)\Delta_J(D_s)$, 
where the tensor product is over $U(\gg)$ (and needs some completion).
Though the noncommutativity of the base and completions
makes it far more complicated,
this is similar to the classical example~(\cite{xu}),
where the algebra $\Diff(M)$ of smooth differential operators on a 
manifold $M$ becomes a bialgebroid over the commutative base $C^\infty(M)$
via the coproduct which multiplicatively extends the 
rule $X\mapsto 1\otimes X + X\otimes 1$ for vector fields $X$ and
$f\mapsto f\otimes 1$ for functions $f$; 
the canonical embedding $C^\infty(M)\hookrightarrow\Diff(M)$
serves both as the source and the target map. 

While we justified the initial formulas in Section~\ref{sec:invdiffop} by 
the formal geometry on a Lie group 
(a wider picture in formal geometry 
will be exhibited in~\cite{SkodaStojicformal}), 
much of the paper is continued in the dual algebraic language 
dictated by the physical motivation
where the Lie algebra generators are interpreted as deformed coordinates, 
rather than invariant vector fields. 
A different variant of this Hopf algebroid structure 
has been outlined in~\cite{tajronkov,tajron,tajronijmpa}, for the special case
when the Lie algebra is the $\kappa$-Minkowski space,
at a physical level of rigor.

\subsection{Algebraic preliminaries}\label{ss:algprel}
We assume familiarity 
with bimodules, coalgebras, comodules, bialgebras, Hopf algebras,
Hopf pairings and the Sweedler notation 
for comultiplications (coproducts) 
$\Delta(h) = \sum h_{(1)}\otimes h_{(2)}$
(with or without the explicit summation sign).
We do not assume previous familiarity with Hopf algebroids. 
In noncommutative geometry, one interprets   
Hopf algebroids~\cite{bohmHbk,Bohm,BrzMilitaru,Lu} 
as formal duals to quantum groupoids. 

The generic symbols for the multiplication map,  
comultiplication, counit and antipode will be $m,\Delta,\epsilon,\mcS$,  
with various subscripts and superscripts. 
All algebras are over a fixed ground field $\genfd$ of characteristic zero 
(in physical applications $\mathbb{R}$ or $\mathbb{C}$). 
The Einstein summation convention on repeated indices is
assumed throughout the article. 
The opposite algebra of an associative algebra $A$ is denoted $A^\op$, 
and the coopposite coalgebra to $C = (C,\Delta)$ is 
$C^{\mathrm{co}} = (C,\Delta^\op)$.
Given a vector space $V$, denote its algebraic
dual by $V^* := \operatorname{Hom}(V,\genfd)$, 
and the corresponding symmetric algebras $S(V)$ and $S(V^*)$. 
If an algebra $A$ is graded, 
we label its graded (homogeneous) 
components by upper indices, $A=\oplus_{i=0}^\infty A^i$, 
$A^i\cdot A^j\subset A^{i+j}$,
and, if $B$ is filtered, we label its filtered 
components $B_0\subset B_1\subset B_2\subset\ldots$ 
by lower indices, $B_i\cdot B_j\subset B_{i+j}$ and $B=\cup_{i=0}^\infty B_i$.
When applied to spaces, we use the hat symbol $\hat{}$ for completions.
We often use the completion of the symmetric algebra 
$\hat{S}(V^*) = \underset\longleftarrow\lim{}_i S_i(V^*)\cong\prod_i S^i(V^*)$ of a Lie algebra $V$ (our main example is when $V$ 
is the underlying space of a Lie algebra $\gg$)
which is the completion of $S(V^*)$ with respect to the degree
of polynomial; it may be identified with the formal power series ring
$\genfd[[\partial^1,\ldots,\partial^n]]$ in $n$ variables.
For our purposes, it is the same to regard this ring, 
as well as the algebraic duals $U(\gg)^*$ and $S(\gg)^*$
of the enveloping and symmetric algebras, 
either as topological or as cofiltered algebras (see Appendix A.2);
the continuous linear maps then translate as linear maps
distributive over formal sums.

The $n$-th {\bf Weyl algebra} $A_n$ 
is the associative algebra generated by 
$x_1,\ldots,x_n$, $\partial^1,\ldots,\partial^n$ subject to relations
$[x_\alpha,x_\beta] = [\partial^\alpha,\partial^\beta] = 0$ and
$[\partial^\alpha,x_\beta] = \delta^\alpha_\beta$. It has a vector
space basis formed by all expressions of the form 
$x_{\alpha_1}\cdots x_{\alpha_k} \partial^{\beta^1}\cdots\partial^{\beta^l}$; 
if we define the degree of this element as $\beta^1+\ldots+\beta^l$,
then $A_n$ becomes a filtered algebra; 
it has no zero divisors and the elements of
the degree at least $k$ form an ideal $(A_n)_{\mathrm{deg}\geq k}$. 
Thus we can form the (semi)completed Weyl algebra 
$\hat{A}_n=\lim_s A_n/(A_n)_{\mathrm{deg}\geq s}$ (``completed by the degree'').
In the geometric part of the paper we shall also consider
the $n$-th {\bf covariant Weyl algebra} $A_n^{\mathrm{cov}}$ 
where the position of the upper versus lower indices in the notation 
will be interchanged; hence $[\partial_\alpha,x^\beta] = \delta_\alpha^\beta$.
Here we shall similarly dually complete by the dual degree which is
$\alpha_1+\ldots+\alpha_k$ on the basis elements 
$x^{\alpha_1}\cdots x^{\alpha_k}\partial_{\beta_1}\cdots\partial_{\beta_l}$ 
to obtain the completion $\hat{A}_n^{\mathrm{cov}}$. The correspondence
$x^\alpha\mapsto \partial^\alpha$ and $\partial_\beta\mapsto x_\beta$ 
extends to the canonical antiisomorphism 
$\hat{A}_n^{\mathrm{cov}}\to\hat{A}_n$.
The {\bf Fock space} is the faithful representation of $A_n$
on the polynomial algebra in $x_1,\ldots,x_n$ where
each $x_\mu$ acts as the multiplication operator and $\partial^\mu$ as 
the partial derivative; this action extends continuously to
a unique action of $\hat{A}_n$ also called Fock.

\begin{defn}
Let $A$ be an algebra and $B$ a bialgebra. 

A left action $\triangleright : B\otimes A\to A$ 
(right action $\triangleleft : A\otimes B\to A$), 
is a left (right) {\bf Hopf action} if 
$b\triangleright(a a') = \sum (b_{(1)}\triangleright a)
(b_{(2)}\triangleright a')$ and $b\triangleright 1 =\epsilon(b)1$
(or, respectively, 
$(aa')\triangleleft b = \sum (a\triangleleft b_{(1)})(a\triangleleft b_{(2)})$
and $1\triangleleft b = \epsilon(b)1$),
for all $a,a'\in A$ and $b,b'\in B$.
We then also say that $A$ is a left (right) $B$-module algebra.
As usual, we freely exchange actions and representations; thus
by abuse of language we say that 
a representation $\psi:B\to\End(A)$ is a left Hopf action (representation)
if $b\otimes a\mapsto \psi(b)(a)$ is a left Hopf action.
Given a left Hopf action, the {\bf smash product} $A\sharp B$ 
(for a right Hopf action, the smash product $B\sharp A$) 
is an associative algebra which is 
a tensor product vector space $A\otimes B$ ($B\otimes A$)
with the multiplication bilinearly extending the formulas
$$\begin{array}{l}
(a\sharp b) (a'\sharp b') = 
\sum a (b_{(1)}\triangleright a') \sharp b_{(2)} b',
\,\,\,\,a,a'\in A, b,b'\in B,
\\
(b\sharp a)(b'\sharp a') = 
\sum b b'_{(1)}\sharp (a\triangleleft b'_{(2)}) a',
\,\,\,\,a,a'\in A, b,b'\in B,
\end{array}
$$
where, for emphasis, one writes $a\sharp b := a\otimes b$.
\end{defn}
Note that $A\sharp 1$ and $1\sharp B$ are subalgebras in $A\sharp B$, 
canonically isomorphic to $A$ and $B$. 
If $B$ is a Hopf algebra with an antipode $\mcS$, we 
may replace a left Hopf action $\psi:B\to\operatorname{End}A$
by a homomorphism
$\psi\circ\mcS:B^{\mathrm{co}}\to\operatorname{End}^{\mathrm{op}}A$,
yielding a {\it right} Hopf action 
$\triangleleft:A\otimes B^{\mathrm{co}}\to A$, 
$\triangleleft:a\otimes b\mapsto 
a\triangleleft b := \mcS(b)\triangleright a$,
and enabling us to define the smash product $B^{\mathrm{co}}\sharp A$.
If $B$ is cocommutative (for instance, $B^{\mathrm{co}}=B=U(\gg)$ below)
then $\mcS^2 = \id$ 
and there is an isomorphism $A\sharp B\cong B\sharp A$ of algebras, 
$a\sharp b\mapsto\sum b_{(1)}\sharp (a\triangleleft b_{(2)})$, 
with the inverse 
$b\sharp a\mapsto \sum (b_{(1)}\triangleright a)\sharp b_{(2)}$. 

Simple examples of smash products are the Weyl algebras $A_n$ 
(and completions $\hat{A}_n$). 
Indeed the symmetric
algebra $S(V)$ of a vector space is a Hopf algebra with $\Delta(x) = 1\otimes x+ x\otimes 1$ for generators $x\in V$; and if $V$ is a vector
space spanned by $x_1,\ldots,x_n$ then there are canonical isomorphisms 
$\hat{A}_n = S(V)\sharp S(V)^*\cong S(V)^*\sharp S(V)$
where the smash products are constructed using the 
(right and left) Hopf actions of $S(V)$ on $S(V)^*$
defined using duality. More generally, replacing $S(V)$ by
its noncommutative generalization -- the universal enveloping algebra
$U(\gg)$ -- we explictly construct in Section~\ref{sec:defphasesp} 
certain smash products $H^L = U(\gg^L)\sharp S(\gg)^*$ and 
$H^R = S(\gg)^*\sharp U(\gg^R)$, both isomorphic as
algebras to $\hat{A}_n$; their special smash product structures however
give rise to a left and a right $U(\gg)$-bialgebroid structures 
(in a completed sense). 

\begin{defn} \cite{bohmHbk,BrzMilitaru}\label{def:leftbialg}
A {\bf left bialgebroid} $(H,m,\alpha,\beta,\Delta,\epsilon)$ 
over the {\bf base algebra} $\AA$ (shortly, left $\AA$-bialgebroid) 
consists of 
\begin{itemize}
\item ({\bf total algebra})
an associative algebra $H$ with multiplication $m$ 
\item ($\AA$-bimodule structure on $H$)
morphisms of algebras
{\bf source} $\alpha:\AA\to H$ and {\bf target} $\beta:\AA^\op\to H$ 
satisfying $[\alpha(a),\beta(b)] = 0$ for all $a,b\in A$, 
hence equipping $H$ with the 
structure of an $\AA$-bimodule via the formula $a.h.b := 
\alpha(a)\beta(b) h$ for $a,b\in \AA$ and $h\in H$;
\item ($\AA$-coring structure on $H$, see Definition~\ref{def:coring}) 
$\AA$-bimodule maps coproduct $\Delta: H\to H\otimes_\AA H$ 
and the corresponding counit $\epsilon : H\to \AA$ 
making $(H,\Delta,\epsilon)$ into a comonoid (coalgebra) 
in the category of $\AA$-bimodules equipped
with the tensor product $\otimes_\AA$ of $\AA$-bimodules.
\end{itemize}
In addition, $\Delta$ and $\epsilon$ need to be compatible with the 
multiplication $m$, but in more subtle way than in the bialgebra case. 
Namely, if the base $\AA$ is noncommutative,
$H\otimes_\AA H$ does not inherit a well defined multiplication
from the usual tensor product $H\otimes H$ over the ground field. 
Instead, one demands that the image of $\Delta$ 
is within the subspace $H\times_\AA H$ 
consisting of all $\sum_i b_i\otimes b'_i$
in $H\otimes_\AA H$ such that $\sum_i b_i\otimes b'_i \alpha(a)
= \sum_i b_i\beta(a)\otimes b'_i$ for all $a\in\AA$;
it appears that $H\times_\AA H$ automatically inherits
the well-defined algebra structure from $H\otimes H$.
We demand that after corestricting $\Delta$
to this smaller codomain $H\times_\AA H$, 
$\Delta$ becomes an algebra map. Similarly, 
the counit $\epsilon$ is not required to be an algebra map, 
but a weaker condition is assumed:
the formula $h\otimes a\mapsto \epsilon(h\alpha(a))$ needs to define 
an action $H\otimes A\to A$ restricting to the multiplication
$A\otimes A\to A$. 
\end{defn}

\subsection{Preliminaries on formal differential operators}\label{ss:difprel}
If $\genfd$ is $\mathbb{R}$ or $\mathbb{C}$ 
then for any smooth  manifold $M$ of dimension $n$
we denote by $C^\infty(M)$ the algebra of smooth functions. 
Given a point $e\in M$ and a natural number $s$, recall that an $s$-jet
of functions around $e$ is a class of equivalence of 
smooth functions defined locally 
around $e$, where two functions are equivalent
if they are defined in some neighborhood of $e$ 
and their Taylor series up to order $s$ agree at $e$. 
All $s$-jets around $e$ form a vector space
$J^s(M,e)$ with canonical projections $J^{s+1}(M,e)\to J^s(M,e)$
and the inverse limit 
$J^\infty(M,e)= \underset{\longleftarrow}\lim{}_{s}\, J^s(M,e)$ 
is by definition the space of formal functions around $e$;
the spaces $J^s(M,e)$ with their canonical projections 
then form a cofiltration of $J^\infty(M,e)$. 
In any chosen coordinate chart around $e$ 
the formal functions are represented 
by formal power series in $n = \mathrm{dim}\,M$ indeterminates. 
Similarly, one can consider $s$-jets of maps to
other manifolds (including the coordinate charts 
viewed as maps to $\genfd^n$) and, in the limit, 
formal maps and formal charts. 
Since, by a theorem of \'E.~Borel~\cite{EBorel}, 
each formal power series over $\mathbb{R}$ 
is a Taylor series of
a non-unique smooth function, a formal function may be viewed as
an $\infty$-jet of an actual but non-unique smooth function. 
Thus those quantities in differential geometry which depend 
only on their Taylor series have formal analogues, 
namely the $\infty$-jets of actual locally defined smooth quantities. 

A regular differential operator $Q\in\Diff_s(M)$ of degree up to $s$ 
is in every smooth chart a sum of the form
$\sum_{|J|\leq s} q^J\partial_J$ where the sum is over multiindices
$J=(j_1,\ldots,j_n)\in\mathbb{N}^n_0$ with $|J| = j_1+\ldots+j_n\leq s$,
and $q^J$ is a smooth function defined over the chart.
At the jet level, the ring of regular differential operators 
$\Diff(M) = \cup_{s\in\mathbb{N}}\Diff_s(M)\subset\End_{\mathbb{R}}(C^\infty(M))$
gives rise to the ring 
$\Diff^\omega(M,e)\subset\End_{\mathbb{R}}(J^\infty(M,e))$
of formal differential operators at $e$, namely the 
$\infty$-jets of regular differential operators around $e$. 
A formal differential operator at $e$ is a sum 
$\sum_{|J|\leq s} q^J \partial_J$ where $q^J = q^J(x^1,\ldots,x^n)$
is a formal function at $e$;  
these sums can be viewed
as elements of the semicompleted Weyl algebra $\hat{A}_n^{\mathrm{cov}}$.
The evaluation of a differential operator at a function at $e$ 
is a rule for a degenerate pairing between
$\Diff^\omega(M,e)$ and $J^\infty(M,e)$. If $M = G$ is a Lie group and
$e\in G$ the unit element, then it 
restricts to a nondegenerate pairing between the subspace 
$\Diff^{\omega,R}(G,e)\subset\Diff^\omega(G,e)$ of right invariant
formal differential operators and $J^\infty(G,e)$. 

For the later transition to the noncommutative point of view, 
it is useful to consider also the algebra of differential operators 
acting to the {\it left},
which is simply the opposite algebra
$\Diff^\op(M)\subset\End_{\mathbb{R}}^\op(C^\infty(M))$
and its formal version $\Diff^{\omega,\op}(M,e)$. In order to stick to
the Weyl algebra notation and commutation relations, 
after changing the order of operators we 
denote $x_\mu = (\partial_\mu)^\op$ and $\partial^\nu = (x^\nu)^\op$. The
canonical antiisomorphism $\Diff^\omega(M,e)\to\Diff^{\omega,\op}(M,e)$
hence sends $\sum_{|J|\leq s} q^J \partial_J$ to 
$\sum_{|J|\leq s} x_J p^J(\partial^1,\ldots,\partial^n)$ 
where $p^J$ is $(q^J)^\op$ (written as a formal function 
of $\partial^\mu$). The latter sum 
can be viewed as belonging to the semicompleted
Weyl algebra $\hat{A}_n$ with contravariant notation 
(as in~\ref{ss:algprel}).

\section{Left versus right invariant differential operators}
\label{sec:invdiffop}

In Section~\ref{sec:defphasesp} we shall introduce 
a noncommutative phase space $H^L$ 
of Lie type and important matrix $\OO$ which plays the central
role in defining our Hopf algebroid structure. 
Geometrical origin of this matrix and related issues 
are clarified in this section using calculations relating
left and right invariant vector fields.

Throughout the article, $\gg$ is a fixed
Lie algebra over $\genfd$ of some finite dimension $n$.
In a basis $\hx_1,\ldots,\hx_n$ of $\gg$, we define the structure 
constants $C_{\mu\nu}^\lambda$ by
\begin{equation}\label{eq:cx}
[\hx_\mu,\hx_\nu] = C_{\mu\nu}^\lambda \hx_\lambda.
\end{equation}
Introduce the opposite Lie algebra $\gg^R$ generated by $\hy_\mu$, where
\begin{equation}\label{eq:cy}
[\hy_\mu,\hy_\nu] = -C_{\mu\nu}^\lambda \hy_\lambda.
\end{equation}
The Lie algebra $\gg^R$ is antiisomorphic to $\gg^L:=\gg$ 
via $\hy_\mu\mapsto \hx_\mu$,
inducing an isomorphism $U(\gg^L)^\op\cong U(\gg^R)$.

If $\genfd$ is $\mathbb{R}$ or $\mathbb{C}$ we also fix a Lie
group $G$ with unit $e$ such that $\gg$ is its Lie algebra, realized
as the algebra $\Vect^L(G)$ of left invariant vector fields on $G$,
then $\gg^R\cong\Vect^R(G)$.
The universal enveloping algebra $\gg^R\hookrightarrow U(\gg^R)$
can be realized as the algebra of 
right invariant differential operators on $G$, i.e.\ by 
embedding $\gg^R\cong\Vect^R(G)\hookrightarrow\Diff^R(G)$. 
If $R_g:G\to G$ is the right multiplication by $g\in G$ 
then a differential operator $D\in\Diff(G)$ is right invariant 
if $(R_{g*})_{h} D_h = D_{h g}$. Therefore $D_g = (R_{g*})_{e}D_e$ and
every right invariant formal differential operator $D\in\Diff^{\omega,R}(G,e)$ 
at the unit $e$ (cf.~\ref{ss:difprel}) 
extends to a unique right invariant analytic 
differential operator on $G$. 
Thus $U(\gg^R)\cong\Diff^{\omega,R}(G,e)\hookrightarrow\Diff^\omega(G,e)$
and the evaluation of differential operator at $\infty$-jets of smooth functions gives a pairing of $U(\gg)$ and $J^\infty(G,e)$ which
is nondegenerate by the L.~Schwartz's theorem~(\cite{CartierHopf}). 

The generators of the universal enveloping algebra $U(\gg^L)$ and $U(\gg^R)$ 
are also denoted by $\hx_\mu$ and $\hy_\mu$, unlike the generators 
of the symmetric algebra $S(\gg)$ which are 
denoted by $x_1,\ldots,x_n$ (without hat symbol) instead. 
Each element $Y\in\gg^R$ can be written as $Y = Y^\mu\hy_\mu$, thus
$Y^\mu:\gg\to\genfd$ may be taken as global coordinates on $\gg^R$ 
as a manifold.
The exponential map $\exp: T_e G\to G$ restricts to
a diffeomorphism from some star-shaped open neighborhood 
$U$ of $0\in\gg^R$ 
to some open neighborhood $V = \exp(U)$ of $e\in G$. 
Thus the functions $w^\mu:V\to\mathbb{R}$ 
together with the basis
$\partial_\nu^w$ of vector fields on $V$ 
given by 
\begin{equation}\label{eq:wcoord}
w^\mu(\exp(Y^\gamma\hy_\gamma))=Y^\mu,\,\,
\partial_\nu^w = (d\exp)(\partial/\partial Y^\nu),\,\,\,\,\,\,
\mu,\nu=1,\ldots,n,
\end{equation}
form a system of coordinates on the tangent manifold $T V$. 
The corresponding multiplication by a coordinate and the 
derivative elements in $\Diff(V)$ satisfy  
the usual commutation relations $[\partial^w_\mu,w^\nu] = \delta^\nu_\mu$ 
generating a copy ${}^w A_n^{\mathrm{cov}}$ 
of the Weyl algebra $A_n^{\mathrm{cov}}$ (see~\ref{ss:algprel}). 
There is a well known 
formula~(see~\cite{helgason}, Chapter II or~\cite{postnikovSemV}, Lecture 4 Cor. 1) 
for the differential 
$d\exp : TU\to TG$ of the exponential map $\exp|_U:U\to G$,
\begin{equation}\label{eq:diffexp}
(d \exp)_Y
= (R_{\exp Y *})_{e}\circ\frac{1-e^{-\ad Y}}{\ad Y}\,\,\,
\mbox{for}\,\,\mbox{all}\,\,\,
Y = Y^\gamma\hy_\gamma\in\gg^R,
\end{equation}
where the action of $\ad Y =\ad^R Y$ is understood in the sense of the
identification $T_Y U\cong \gg^R$ of the tangent space at $Y$ with $\gg^R$
(hence it is $-\ad Y$ in the sense of $\gg^L$-bracket).
Let ${}^w\CC$ be the matrix of functions
$({}^w\CC)^\mu_\nu = C^\mu_{\nu\gamma} w^\gamma = - C^\mu_{\gamma\nu}w^\gamma:V\to\mathbb{R}$.
For fixed $Y = Y^\mu\hy_\mu\in\gg^R$, the calculation
$(\ad Y)(\hy_\nu) = [Y^\mu\hy_\mu,\hy_\nu]_{\gg^R} = - Y^\mu C^\gamma_{\mu\nu}\hy_\gamma$ 
implies
$$
(\ad Y)^N(\hy_\beta) = ({}^w\CC^N(\exp{Y}))^\gamma_\beta\hy_\gamma,
\,\,\,\,\,w^\mu(\exp{Y})=Y^\mu,\,\,\,N= 0,1,2,\ldots
$$ 
By~(\ref{eq:diffexp}) we have
$(R_{\exp{Y}*})_e \hy_\alpha = 
(d \exp)_Y\circ\frac{\ad Y}{1-\exp(-\ad Y)}\hy_\alpha$,
which equals
$(d\exp)_Y \left(\frac{-{}^w\CC}{e^{-{}^w\CC}-1}\right)^\beta_\alpha\hy_\beta
= \left(\frac{-{}^w\CC}{e^{-{}^w\CC}-1}\right)^\beta_\alpha(d\exp)_Y\hy_\beta$;
hence the basis $\hy_\alpha:w\mapsto (R_{\exp{Y}*})_{e} (\hy_\alpha)$
of the space of right invariant vector fields $\Vect^R(G)|_V$ is in 
the coordinates $w^1,\ldots,w^n,\partial^w_1,\ldots,\partial^w_n$ given by 
\begin{equation}\label{eq:univwC}
\hy^{\exp}_\alpha = (R_{\exp{Y}*})_e(\hy_\alpha) = 
\left(\frac{-{}^w\CC}{e^{-{}^w\CC}-1}\right)^\beta_\alpha \partial_\beta^w.
\end{equation}
Notice that $\frac{-{}^w\CC}{e^{-{}^w\CC}-1}$ is a matrix of power series
in $w^1,\ldots,w^n$, hence analytic in $V$.
The map $()^{\exp}:\hy_\mu\mapsto \hy^{\exp}_\mu$ 
is an embedding of $U(\gg)$ into the 
algebra of formal differential operators $\Diff^\omega(U,e)$ 
with the distinguished Weyl subalgebra ${}^w A_n^{\mathrm{cov}}$
in the coordinate chart $w^1,\ldots,w^n$. 
The same geometric embedding is obtained by Durov 
in a more general setting of formal
geometry over general ring $\genfd\supset\mathbb{Q}$ in~\cite{ldWeyl}, 
formula (36), where ${}^w\CC$ is denoted by $M$. 
Notice that $L_{\exp{Y}} = L_{\exp{Y}}R_{\exp{(-Y)}}R_{\exp{Y}} = 
R_{\exp{Y}}L_{\exp{Y}}R_{\exp{Y}}^{-1}$, hence using $\Ad_{\exp{Y}} = e^{\ad{Y}}$
we obtain
$$
(L_{\exp{Y}*})_e = {}^w\OO(\exp{Y}) \circ (R_{\exp{Y*}})_e
= (R_{\exp{Y}*})_e\circ\Ad_{\exp{Y}},
$$
$$
(L_{\exp{Y}*})_e(\hy_\alpha) = ((R_{\exp{Y}*})_e\circ e^{\ad{Y}})(\hy_\alpha),
$$
where
\begin{eqnarray}
({}^w\OO)(\exp{Y}) 
= (L_{\exp{Y}*})_e\circ (R_{\exp{(-Y)}*})_{\exp{Y}}:T_{\exp{Y}}G\to T_{\exp{Y}}G,
\\
\Ad_{g} 
:= (L_{g*})_{g^{-1}} \circ(R^{-1}_{g*})_e:\gg^R\to\gg^R,\,\,\,\,\,\,g\in G,
\\
e^{{}^w\CC}(\exp Y) = e^{\ad{Y}} = \Ad_{\exp{Y}},
\\
{}^w\OO^\beta_\alpha := (e^{{}^w\CC})^\beta_\alpha,
\end{eqnarray}
hence the bases of $\Vect^R(G)|_V$ and $\Vect^L(G)|_V$ 
are related via ${}^w\OO^\beta_\alpha$,
\begin{equation}\label{eq:hxOhy}
\hx_\alpha^{\exp{}} 
= {}^w\OO^\beta_\alpha \hy_\alpha^{\exp{}}
\end{equation}
\begin{equation}\label{eq:univwCx}
\hx^{\exp}_\alpha := (L_{\exp{Y}*})_e(\hy_\alpha) = 
\left(\frac{{}^w\CC}{e^{{}^w\CC}-1}\right)^\beta_\alpha \partial_\beta^w.
\end{equation}

\section{From differential operators to 
the deformed phase space}
\label{sec:defphasesp}

If $C^\lambda_{\mu\nu}=0$ then 
$\hx_\mu^{\exp{}} = \hy_\mu^{\exp{}} = \partial_\mu^w$, 
which is not in the spirit of the interpretation in physics 
where $\hx_\mu$ are often viewed as the analogue or deformation
of commutative coordinates $x_\mu$, 
cf.~\cite{AC,Govstat,HallidaySzabo,tajronkov,tajron}. 
For that purpose most of the paper is written in 
somewhat dual language obtained as follows.
Introduce  the {\it antiisomorphism}
$\Diff^\omega(G,e)\to\hat{A}_n$ 
(restricting to ${}^w A_n^{\mathrm{cov}}\to A_n$)
mapping $w^\mu\mapsto\partial^\mu$ and $\partial_\nu^w\mapsto x_\nu$ and 
consequently ${}^w\CC\mapsto\CC$, 
$\frac{-{}^w\CC}{e^{-{}^w\CC}-1}\mapsto\phi$, 
$\hy^{\exp}_\alpha\mapsto\hy^\phi_\alpha := x_\beta\phi^\beta_\alpha$, 
$\frac{{}^w\CC}{e^{{}^w\CC}-1}\mapsto\tilde\phi$,
$\hx^{\exp{}}\mapsto \hx^\phi := x_\beta\tilde\phi^\beta_\alpha$, 
${}^w\OO\mapsto\OO$,
where $\CC$, $\OO$, $\phi$ and $\tilde\phi$  are $n\times n$ matrices
\begin{equation}\label{CCmatrix}
\CC^\alpha_\beta := C^\alpha_{\beta\gamma}\partial^\gamma,
\,\,\,\,\,\,\alpha,\beta=1,\ldots,n,
\,\,\,\,\,\,\,\OO := e^\CC,
\end{equation}
\begin{equation}\label{eq:phi}
\phi := \frac{-\CC}{e^{-\CC} -1} = 
\sum_{N=0}^\infty \frac{(-1)^N B_N}{N!} \CC^N,\,\,\,\,\,\,\,\,\,\,\,\,\,\,
\tilde\phi = \frac{\CC}{e^\CC-1}.
\end{equation}
The constants $B_N$ are the Bernoulli numbers and the matrix entries 
$\phi^\beta_\alpha$, $\tilde\phi^\beta_\alpha,\OO^\mu_\nu\in\hat{S}(\gg^*)$ 
are formal power series in the elements $\partial^1,\ldots,\partial^n$
of $S(\gg)^*$ 
which correspond to the basis of $\gg^*$ dual to 
$\hx_1,\ldots,\hx_n$ of $\gg^L$. This is in agreement with
the notation in the Weyl algebra $A_n$ and in
$\hat{A}_n\cong S(V)\sharp S(V)^*$ for $V=\gg$. 
The formula $\bm\phi_+(\hat{x}_\alpha)(\partial^\beta) := \phi^\beta_\alpha$
determines a linear map $\bm\phi_+(\hat{x}_\alpha):\gg^*\to\hat{S}(\gg^*)$,
which by the Leibniz rule and continuity extends to a unique 
continuous derivation
$\bm\phi_+(\hx_\alpha)\in\operatorname{Der}(\hat{S}(\gg^*))$. 
It is crucial that $\hx_\alpha\mapsto\bm\phi_+(\hx_\alpha)$ 
defines a Lie algebra homomorphism
$\bm\phi_+:\gg^L\to\operatorname{Der}(\hat{S}(\gg^*))$.
Equivalently, $\bm\phi_+$ 
extends to a unique right {\it Hopf} action also denoted 
\begin{equation}\label{eq:bmphi}
\bm\phi_+: U(\gg^L)\to\operatorname{End}^{\mathrm{op}}(\hat{S}(\gg^*)).
\end{equation}
This induces the smash product 
$H^L := U(\gg^L)\sharp_{\bm\phi_+}\hat{S}(\gg^*)$
interpreted as the 'noncommutative phase space of Lie type'.
(Warning: in~\cite{heisd} we used the notation $\bm\phi$ 
for the {\it left} Hopf action $\bm\phi_- = \bm\phi_{+}\circ\mcS_{U(\gg^L)}$,
where $\mcS_{U(\gg^L)} = \mcS^{-1}_{U(\gg^L)}$ is the antipode for
$U(\gg^L)$, satisfying $\gg^L\ni h\mapsto -h$).

Regarding that $\gg^R$ is a Lie algebra with known structure constants, 
$-C^\alpha_{\beta\gamma}$, the formula~(\ref{eq:phi}) can be applied
to it. This also gives the right Hopf action $\tilde{\bm\phi}_+:U(\gg^R)
\to\operatorname{End}^\op(\hat{S}(\gg^*))$, $\tilde{\bm\phi}_+(\hy_\nu)
(\partial^\mu)=\tilde\phi^\mu_\nu$; the right bialgebroid structure
constructed below will however be based on the {\it left} Hopf
action $\tilde{\bm\phi}_- = \tilde{\bm\phi}_+\circ\mcS_{U(\gg^R)}:U(\gg^R)
\to\operatorname{End}(\hat{S}(\gg^*))$, $\tilde{\bm\phi}_-(-\hy_\nu)
(\partial^\mu)=\tilde\phi^\mu_\nu$. Thus we can define the
smash product $H^R := \hat{S}(\gg^*)\sharp_{\tilde{\bm\phi}_-}U(\gg^R)$.
Its generators 
are $\hat{y}_\mu,\partial^\mu$, $\mu=1,\ldots,n,$ completing
in $\partial^\mu$-s.
In addition to the relations in $U(\gg^R)$ and $\hat{S}(\gg^*)$,
we also have
$$
[\partial^\mu,\hy_\nu] = \left(\frac{\CC}{e^{\CC}-1}\right)^\mu_\nu.
$$

Precomposing $()^{\exp{}}:U(\gg^R)\to\Diff^\omega(G,e)$
by the antiisomorphism $U(\gg^L)\to U(\gg^R)$, $\hx_\mu\mapsto\hy_\mu$
and postcomposing by the above antiisomorphism
$\Diff^\omega(G,e)\to\hat{A}_n$ we obtain 
the monomorphism $U(\gg^L)\to U(\gg^R)\to\Diff^{\omega,R}(G,e)\to\hat{A}_n$
denoted $()^\phi:U(\gg^L)\to\hat{A}_n$, used in 
the rest of the article and called the  
$\phi$-{\bf realization} of $U(\gg^L)$
(by dually-formal differential operators).
When complemented by the rule $\partial^\mu\mapsto\partial^\mu$, 
the $\phi$-realization extends to
a unique continuous isomorphism of algebras 
$U(\gg^L)\sharp_{\bm\phi_+}\hat{S}(\gg^*)\cong \hat{A}_n$,
the $\phi$-realization of $H^L$. 
Notice that $(\hat{x}_\nu)^\phi = \hat{x}_\nu^\phi = x_\rho\phi^\rho_\nu$.
We commonly identify $\hat{S}(\gg^*)$ 
with the subalgebra $1\sharp\hat{S}(\gg^*)$
and $U(\gg^L)$ with $U(\gg^L)\sharp 1$. It follows that in $H^L$
\begin{equation}\label{eq:parhx}
[\partial^\mu,\hx_\nu] = \left(\frac{-\CC}{e^{-\CC} -1}\right)^\mu_\nu.
\end{equation}
This identity justifies the interpretation of $\partial^\mu$ within $H^L$ as
deformed partial derivatives. 
The universal formula~(\ref{eq:phi}) for $\phi$ is, 
in this context, derived in~\cite{ldWeyl} and $H^L$ 
is studied in~\cite{scopr}. 

The map $J^\infty(G,e)\to\hat{S}(\gg^*)$, $w^\nu\mapsto\partial^\nu$
is an antiisomorphism of algebras and it can be combined with the
realization of $U(\gg^L)^\op$ via $\Vect^R(G,e)$ 
to compare with the opposite smash product algebra,  
$${}^w\hat{A}_n^{\mathrm{cov}}\cong\Diff^\omega(G,e)
\cong J^\infty(G,e)\sharp\Vect^R(G,e)\cong
(U(\gg^L)\sharp_{\bm\phi_+}\hat{S}(\gg^*))^\op.$$
The smash product $J^\infty(G,e)\sharp\Vect^R(G,e)$ could be also directly
observed using the duality between $J^\infty(G,e)$ and $\Vect^R(G,e)$.

Similarly to the $\phi$-realization of $U(\gg^L)$, 
there is a $\tilde\phi$-realization of $U(\gg^R)$
exteding to an isomorphism $H^R\cong\hat{A}_n$
given by $\hy_\nu\mapsto x_\rho\tilde\phi^\rho_\nu$, 
$\hat{S}(\gg^{R*})\ni\partial^\nu\mapsto \partial^\nu\in\hat{A}_n$.



\begin{theorem} \label{thO}
There is a unique algebra isomorphism from 
$H^L=U(\gg^L)\sharp\hat{S}(\gg^*)$ 
to $H^R=\hat{S}(\gg^*)\sharp U(\gg^R)$ 
which fixes the commutative subalgebra $\hat{S}(\gg^*)$ 
(i.e. identifies $1\sharp\hat{S}(\gg^{L*})$ with 
$\hat{S}(\gg^{R*})\sharp 1$, 
$1\sharp\partial^\mu\mapsto\partial^\mu\sharp 1$), 
and which maps $\hx_\nu \mapsto\hy_\sigma\OO^\sigma_\nu$,
where $\OO = e^{\CC}$ is an invertible
$n\times n$-matrix with entries $\OO^\mu_\nu\in\hat{S}(\gg^*)$ 
and inverse $\OO^{-1} = e^{-\CC}$. 
After the identification, $[\hx_\mu,\hy_\nu]=0$. 
Consequently, the images of $U(\gg^L)\hookrightarrow H^L$ 
and  $U(\gg^R)\hookrightarrow H^R$ mutually commute. 
The following identities hold
\begin{equation}\label{eq:Oy}
[\OO^\lambda_\mu,\hy_\nu]= C_{\rho\nu}^\lambda\OO^\rho_\mu
\end{equation}
\begin{equation}\label{eq:Ox}
[\OO^\lambda_\mu,\hx_\nu] = C_{\mu\nu}^\rho\OO^\lambda_\rho
\end{equation}
\begin{equation}\label{eq:Omx}
[(\OO^{-1})^\lambda_\mu,\hx_\nu]= -C^\lambda_{\rho\nu}(\OO^{-1})^\rho_\mu
\end{equation}
\begin{equation}\label{eq:Omy}
[(\OO^{-1})^\lambda_\mu,\hy_\nu]= -C^\rho_{\mu\nu}(\OO^{-1})^\lambda_\rho
\end{equation}
\begin{equation}\label{eq:COOCO}
C_{\mu\nu}^\tau\OO^\lambda_\tau = C^\lambda_{\rho\sigma}\OO^\rho_\mu\OO^\sigma_\nu,
\,\,\,\,\,\,C_{\mu\nu}^\tau(\OO^{-1})^\lambda_\tau = 
C^\lambda_{\rho\sigma}(\OO^{-1})^\rho_\mu(\OO^{-1})^\sigma_\nu.
\end{equation}
\end{theorem}

\begin{pf} 
The isomorphism $H^L\cong H^R$ is the composition of the two isomorphisms, 
supplied by $\phi$- and $\tilde\phi$-realizations 
$H^L\cong\hat{A}_n\cong H^R$.
If we express $\hx_\mu$ and $\hy_\nu$ within 
$\hat{A}_n$ as $x_\rho\phi^\rho_\mu$ 
and $x_\sigma\tilde\phi^\sigma_\nu$ respectively,
the commutation relation $[\hx_\mu,\hy_\nu]=0$
becomes $[x_\rho\phi^\rho_\mu,x_\sigma\tilde\phi^\sigma_\nu]=0$,
which is the Proposition~\ref{prop:xycommute} (Appendix A.1).
If $\genfd=\mathbb{R}$ this also easily follows
using the antiisomorphism with the geometric picture 
in Section~\ref{sec:invdiffop} where $[\hy_\mu^{\exp{}},\hx_\nu^{\exp{}}] = 0$
because the left and right invariant vector fields commute.
Comparing $\phi$ and $\tilde\phi$
(or using~(\ref{eq:hxOhy})), note that
\begin{equation}\label{eq:phitildeC}
\tilde\phi = \phi\, e^{-\CC},\,\,\,\,\,\,\,\,
\hx_\nu = \hy_\mu (e^\CC)^\mu_\nu=\hy_\mu\OO^\mu_\nu.
\end{equation}
Rewrite $[\hx_\mu,\hy_\nu]$ now as
$$
[\hy_\rho\OO^\rho_\mu,\hy_\nu] = [\hy_\rho,\hy_\nu]\OO^\rho_\mu +
\hy_\lambda[\OO^\lambda_\mu,\hy_\nu] = \hy_\lambda(- C_{\rho\nu}^\lambda\OO^\rho_\mu
+[\OO^\lambda_\mu,\hy_\nu]).
$$
Starting with the evident fact $[\partial^\gamma,\hy_\nu] \in \hat{S}(\gg^*)$,
and using the induction, one
shows $[\hat{S}(\gg^*),\hy_\nu]\subset\hat{S}(\gg^*)$. 
Thus, $(- C_{\rho\nu}^\lambda\OO^\rho_\mu
+[\OO^\lambda_\mu,\hy_\nu])\in\hat{S}(\gg^*)$. 
Elements $\hy_\lambda$ are independent 
in $H^R$, which is here considered a right $\hat{S}(\gg^*)$-module, 
hence $0 = \hy_\lambda(- C_{\rho\nu}^\lambda\OO^\rho_\mu
+[\OO^\lambda_\mu,\hy_\nu])$ implies~(\ref{eq:Oy}). 
Similarly, in $[\hx_\mu,\hy_\nu]=0$ replace
$\hy_\nu$ with $\hx_\lambda(\OO^{-1})^\lambda_\nu$ to prove~(\ref{eq:Omx}).
To show~(\ref{eq:Ox}), calculate $C^\lambda_{\mu\nu}\hy_\rho\OO^\rho_\lambda
= C^\lambda_{\mu\nu}\hx_\lambda = [\hx_\mu,\hx_\nu]= 
[\hy_\rho\OO^\rho_\mu,\hx_\nu] = \hy_\rho [\OO^\rho_\mu,\hx_\nu]$, hence
$\hy_\rho(C^\lambda_{\mu\nu}\OO^\rho_\lambda-[\OO^\rho_\mu,\hx_\nu])=0$. 
For~(\ref{eq:Omy}) we reason analogously with 
$[\hx_\rho(\OO^{-1})^\rho_\mu,\hy_\nu]$. If in~(\ref{eq:Oy}) 
and~(\ref{eq:Omx}) we replace $\hy_\nu$ (resp. $\hx_\nu$) on the left by
$\hy_\rho(\OO^{-1})^\rho_\nu$ (resp. $\hx_\rho\OO^\rho_\nu$), we get a 
quadratic (in $\OO$ or $\OO^{-1}$) expression on the right,
which is then compared with~(\ref{eq:Ox}) and~(\ref{eq:Omy}) to obtain
~(\ref{eq:COOCO}).
\end{pf}

\section{Actions $\blacktriangleright$ and $\blacktriangleleft$ 
and some identities for them}

There is a map $\epsilon_S : \hat{S}(\gg^*)\to\genfd$, 
taking a formal power series to its constant term 
('evaluation at $0$').
We introduce the 'black action' $\btr$
of $H^L$ on $U(\gg^L)$ as the composition
\begin{equation}\label{eq:blact}
H^L\otimes U(\gg^L)\hookrightarrow H^L\otimes H^L
\stackrel{m}\longrightarrow H^L \cong U(\gg^L)\sharp_{\bm\phi_+}\hat{S}(\gg^*)
\stackrel{\id\sharp\epsilon_S}\longrightarrow U(\gg^L),
\end{equation}
where $m$ is the multiplication map. 
$\btr$ is the unique action
for which $\partial^\mu\btr 1 = 0$ for all $\mu$ 
and $\hat{f}\btr 1 = \hat{f}$ for all $\hat{f}\in U(\gg^L)$. 
It follows that $\OO^\mu_\nu\btr 1=\delta^\mu_\nu 1 =(\OO^{-1})^\mu_\nu\btr 1$
and  $\hy_\nu\btr 1 = \hx_\mu (\OO^{-1})^\mu_\nu\btr 1 = 
\delta^\mu_\nu \hx_\mu = \hx_\nu$.
Similarly, the right black action 
$\btl$ of $H^R$ on $U(\gg^R)$ is the composition
$$
U(\gg^R)\otimes H^R \hookrightarrow H^R\otimes H^R 
\stackrel{m}\longrightarrow H^R\cong \hat{S}(\gg^*)\sharp U(\gg^R)
\stackrel{\epsilon_S\sharp \id}\longrightarrow U(\gg^R),
$$ 
characterized by $1\btl\partial^\mu = 0$, 
and $1\btl\hat{u} = \hat{u}$, for all $\hat{u}\in U(\gg^R)$.
The actions $\btr,\btl$ and the smash products $H^L,H^R$ can be described
abstractly 
in terms of the pairings between $\hat{S}(\gg^*)$ and $U(\gg^L)$
or $U(\gg^R)$,
or equivalently in the geometric picture, between $J^\infty(G,e)$ and $\Vect^L(G)$ or
$\Vect^R(G)$~(\cite{StojicPhD}), 
but we stay here within a more explicit approach. 
\begin{theorem}\label{tm:Oind} For any $\hat{f},\hat{g}\in U(\gg^L)$
the following identities hold
\begin{equation}\label{eq:leibnizOx}
\hx_\alpha \hat{f} =
(\OO^\beta_\alpha\btr\hat{f})\hx_\beta
\end{equation}
\begin{equation}\label{eq:DeltaO}
\OO^\gamma_\alpha\btr (\hat{g}\hat{f}) =
(\OO^\beta_\alpha\btr\hat{g})
(\OO^\gamma_\beta\btr\hat{f})
\end{equation}
\begin{equation}\label{eq:DeltamO}
(\OO^{-1})^\gamma_\alpha\btr (\hat{g}\hat{f}) =
((\OO^{-1})^\gamma_\beta\btr\hat{g})
((\OO^{-1})^\beta_\alpha\btr\hat{f})
\end{equation}
\begin{equation}\label{eq:y11x}
\hy_\alpha\btr \hat{f} = \hat{f}\hx_\alpha
\end{equation}
\begin{equation}\label{eq:leibnizOxFull}
(\hx_\alpha \btr\hat{f})\hat{g} =
(\OO^\beta_\alpha\btr\hat{f})(\hx_\beta\btr\hat{g})
\end{equation}
\end{theorem}

\begin{pf} 
We show (\ref{eq:leibnizOx}) for monomials $\hat{f}$ by induction on
the degree of monomial; by linearity this is sufficient. For the base of
induction, it is sufficient to note
$\OO^\beta_\alpha\btr 1 = \delta^\beta_\alpha$. For the step
of induction, calculate for arbitrary $\hat{f}$ of degree $k$
$$\begin{array}{lcl}
\OO^\gamma_\alpha\btr (\hx_\nu\hat{f}) &=&
[\OO^\gamma_\alpha,\hx_\nu \hat{f}]\btr 1
+ \hx_\nu\hat{f}\OO^\gamma_\alpha\btr 1
\\
&=&[\OO^\gamma_\alpha,\hx_\nu]\btr\hat{f}+\hx_\nu[\OO^\gamma_\alpha,\hat{f}]\btr 1
+\hx_\nu\hat{f}\delta^\gamma_\alpha
\\
&=&C^\beta_{\alpha\nu}\OO^\gamma_\beta\btr \hat{f}
+\hx_\nu(\OO^\gamma_\alpha\btr \hat{f})
\\
&=&(C^\beta_{\alpha\nu}+\delta^\beta_\alpha\hx_\nu)(\OO^\gamma_\beta\btr\hat{f})
\\
&=& (\OO^\beta_\alpha\btr\hx_\nu)(\OO^\gamma_\beta\btr\hat{f}),
\end{array}$$
and use this result in the following:
$$\begin{array}{lcl}
\hx_\alpha\hx_\nu\hat{f}& = &(\OO^\beta_\alpha\btr\hx_\nu)\hx_\beta\hat{f}
\\
&=& (\OO^\beta_\alpha\btr\hx_\nu)(\OO^\gamma_\beta\btr\hat{f})\hx_\gamma
\\
&=& (\OO^\gamma_\alpha\btr(\hx_\nu\hat{f}))\hx_\gamma.
\end{array}$$
Thus (\ref{eq:leibnizOx}) holds for $\hat{f}$-s of
degree $k+1$, hence, by induction, for all. Along the way,
we have also shown (\ref{eq:DeltaO}) for $\hat{g}$ of degree $1$ and
$\hat{f}$ arbitrary. Now we do induction on the degree of $\hat{g}$: replace
$\hat{g}$ with $\hx_\mu\hat{g}$ and calculate
$$\begin{array}{lcl}
\OO^\gamma_\alpha\btr((\hx_\mu\hat{g})\hat{f}) &=&
(\OO^\beta_\alpha\btr\hx_\mu) (\OO^\gamma_\beta\btr (\hat{g}\hat{f}))
\\&=&
(\OO^\beta_\alpha\btr\hx_\mu) (\OO^\sigma_\beta\btr \hat{g})(\OO^\gamma_\sigma\btr
\hat{f})
\\&=&
(\OO^\sigma_\alpha\btr(\hx_\mu\hat{g}))(\OO^\gamma_\sigma\btr\hat{f}).
\end{array}$$

The proof of (\ref{eq:DeltamO}) is similar to (\ref{eq:DeltaO}) 
and left to the reader. 
To show (\ref{eq:y11x}), we use (\ref{eq:leibnizOx}) and
the equality $\hy_\alpha = \hx_\beta(\OO^{-1})^\beta_\alpha$ in $H^L$:
$$\begin{array}{lcl}
\hx_\beta(\OO^{-1})^\beta_\alpha\btr \hat{f} &=&
\hx_\beta\btr((\OO^{-1})^\beta_\alpha\btr \hat{f}) =
(\OO^\gamma_\beta\btr((\OO^{-1})^\beta_\alpha\btr
\hat{f}))\hx_\gamma
\\
&=& ((\OO^\gamma_\beta(\OO^{-1})^\beta_\alpha)\btr\hat{f})\hx_\gamma
= \delta^\gamma_\alpha\hat{f}\hx_\gamma = \hat{f}\hx_\alpha
\end{array}$$

Finally, (\ref{eq:leibnizOxFull}) follows from
(\ref{eq:leibnizOx}) by multiplying from the right with $\hat{g}$,
and using $\hx_\beta\btr\hat{g} =\hx_\beta\hat{g}$ and
$\hx_\alpha\btr\hat{f} = \hx_\alpha\hat{f}$.
\end{pf}

Now we state an analogue of the Theorem~\ref{tm:Oind} for $\btl$.

\begin{theorem} \label{th:btl}
For any $\hat{f},\hat{g}\in U(\gg^R)$
the following identities hold
\begin{equation}\label{eq:leibnizOy}
 \hat{f} \hy_\alpha = \hy_\beta
(\hat{f}\btl(\OO^{-1})^\beta_\alpha),
\end{equation}
\begin{equation}\label{eq:DeltaRO}
(\hat{g}\hat{f})\btl\OO^\gamma_\alpha  =
(\hat{g}\btl\OO^\beta_\alpha) (\hat{f}\btl\OO^\gamma_\beta),
\end{equation}
\begin{equation}\label{eq:DeltaRmO}
 (\hat{g}\hat{f})\btl(\OO^{-1})^\gamma_\alpha =
(\hat{g}\btl(\OO^{-1})^\gamma_\beta)
(\hat{f}\btl(\OO^{-1})^\beta_\alpha)
\end{equation}
\begin{equation}\label{eq:z11y}
\hat{f}\btl\hat{z}_\alpha = \hy_\alpha\hat{f},
\end{equation}
\begin{equation}
\hat{g}(\hat{f}\btl\hy_\alpha ) =
(\hat{g}\btl\hy_\beta)(\hat{f}\btl(\OO^{-1})^\beta_\alpha),
\end{equation}
where 
\begin{equation}\label{eq:defz}
\hat{z}_\alpha := \OO^\beta_\alpha \hy_\beta 
= \OO^\beta_\alpha \hx_\rho(\OO^{-1})^\rho_\beta \in H^L\cong H^R.
\end{equation}
\begin{equation}\label{eq:zz}
[\hat{z}_\alpha,\hat{z}_\beta] = C^\gamma_{\alpha\beta}\hat{z}_\gamma
\end{equation}
\end{theorem}

\section{Completed tensor product and bimodules}
\label{sec:bimodule}
In this section, we discuss the completed tensor products 
needed for the coproducts ($\Delta_{S(\gg^*)}$
in this and $\Delta^L$ and $\Delta^R$ in the next section), 
introduce the maps $\alpha^L,\beta^L,\alpha^R,\beta^R$ and use them 
to define $U(\gg^L)$-bimodule structure on $H^L$ and 
$U(\gg^R)$-bimodule structure on $H^R$. 

Note that $S(\gg)=\oplus_{i=0}^\infty S^i(\gg)=\cup_{i=0}^\infty S_i(\gg)$ 
carries a graded and $U(\gg) = \cup_i U_i(\gg)$ 
a {\it filtered} Hopf algebra structure. Both
structures are induced along quotient maps 
from the tensor bialgebra $T(\gg)$. 
By the PBW theorem, the linear map 
\begin{equation}\label{eq:xi}
\xi:S(\gg)\to U(\gg),\,\,\,\,\,\,\,\,
x_{i_1}\cdots x_{i_r}\mapsto \frac{1}{r!}
\sum_{\sigma\in\Sigma(r)}\hx_{i_{\sigma(1)}}\cdots\hx_{i_{\sigma(r)}},
\end{equation} 
is an isomorphism of filtered coalgebras whose inverse $\xi^{-1}$ may be identified
with the projection to the associated graded ring~\cite{Bourbaki,CartierHopf,ldWeyl}.
The isomorphism $\xi$ is related to the $\phi$-realization from
Section~\ref{sec:defphasesp} (hence to the exponential map in the geometric
picture in Section~\ref{sec:invdiffop}) as follows. Consider the Fock
action $\triangleright$ of $\hat{A}_n$ on $S(\gg)$ 
and the $\phi$-realization $()^\phi:U(\gg)\to\hat{A}_n$. 
For each $f,g\in S(\gg)$, $\xi(f)\cdot_{U(\gg)}\xi(g)
= \xi(\xi(f)^\phi\triangleright g)$ and this property uniquely
characterizes $\xi$.

For a multiindex $K=(k_1,\ldots,k_n)\in\mathbb{N}^n_0$, 
denote $|K|:= k_1+\ldots + k_n$, $x_K := x_1^{k_1}\cdots x_n^{k_n}$ and
$\hat{x}_K := \hx_1^{k_1}\cdots \hx_n^{k_n}$. 
The multiindices add up componentwise. 
The partial order on $\mathbb{N}^n_0$
induced by the componentwise $<$ is also denoted $<$.
If $J,K$ are multiindices
the rule $\langle x_k,\partial^J\rangle := J! \,\delta^J_K$ 
continuously in the first factor and linearly extends to a unique 
map $\langle,\rangle:S(\gg)\otimes\hat{S}(\gg^*)\to\genfd$
which is a nondegenerate pairing,
hence it identifies $\hat{S}(\gg^*)\cong S(\gg)^*$. 
This is the unique Hopf pairing extending the duality between $\gg$
and $\gg^*$ where $\hat{S}(\gg^*)$ is the topological Hopf algebra
with elements in $\gg^*$ primitive. By duality, the linear map 
\begin{equation}\label{eq:xiT}
\xi^T:U(\gg)^*\longrightarrow S(\gg)^*\cong\hat{S}(\gg^*)
\end{equation}
transpose (dual) to $\xi$ (see~(\ref{eq:xi})) 
is an isomorphism of cofiltered algebras.

The inclusions of filtered components  
$U_k(\gg)\subset U_{k+1}(\gg)\subset U(\gg)$ induce epimorphisms
of dual vector spaces $U(\gg)^*\to U_{k+1}(\gg)^*\to U_{k}(\gg)^*$,
hence a complete {\it cofiltration} on 
$U(\gg)^* = \underset{\longleftarrow}\lim{}_k\, U_k(\gg)^*$ 
(see Appendix A.2).
For each finite level $k$, $U_k(\gg)$ is finite dimensional, hence 
$(U_k(\gg)\otimes U_l(\gg))^*\cong U_k(\gg)^*\otimes U_l(\gg)^*$.
Thus the multiplication $U_k(\gg)\otimes U_l(\gg)\to U_{k+l}(\gg)
\subset U(\gg)$ 
dualizes to $\Delta_{k,l}: U(\gg)^*\to U_k(\gg)^*\otimes U_l(\gg)^*$.
The inverse limits $\underset\longleftarrow\lim{}_k\, \Delta_{k,k}$
and $\underset\longleftarrow\lim{}_p\,\underset\longleftarrow\lim{}_q\,\Delta_{p,q}$ agree and define the coproduct
$\Delta_{U(\gg)^*} := 
\underset\longleftarrow\lim{}_k\, \Delta_{k,k}  : U(\gg)^*\to
\underset\longleftarrow\lim{}_k\, U_k(\gg)^*\otimes U_k(\gg)^*
\cong \underset\longleftarrow\lim{}_p\underset\longleftarrow\lim{}_q 
U_p(\gg)^*\otimes U_q(\gg)^*$. 
The right-hand side is by definition the completed tensor product,  
$U(\gg)^*\hat\otimes U(\gg)^*$.
(For completed tensoring of {\it elements} and {\it maps} 
we below often use simplified notation, $\otimes$.)
Coproduct $\Delta_{U(\gg)^*}$ transfers, 
along the isomorphism 
$\xi^T:U(\gg)^*\stackrel\cong\longrightarrow S(\gg)^*$ 
of cofiltered algebras (see~(\ref{eq:xiT})), 
to the topological coproduct on the completed symmetric algebra 
$\hat{S}(\gg^*)\cong S(\gg)^*$ (cf.~\cite{scopr}),
$$
\Delta_{\hat{S}(\gg^*)} : \hat{S}(\gg^*)\to \hat{S}(\gg^*)\hat\otimes
\hat{S}(\gg^*).
$$
This construction can be performed both for $\gg^L$ and $\gg^R$.
The canonical isomorphism of Hopf algebras $U(\gg^R)\cong U(\gg^L)^\op$
induces the isomorphism of dual cofiltered Hopf algebras 
$U(\gg^R)^*\cong (U(\gg^L)^*)^{\mathrm{co}}$, commuting with $\xi^T$, 
hence inducing an isomorphism of Hopf algebras
$\hat{S}(\gg^{R*})\cong\hat{S}(\gg^{L*})^{\mathrm{co}}$
fixing the underlying algebra $\hat{S}(\gg^*)$. Thus, the
coproduct on $\hat{S}(\gg^{R*})$ is $\Delta_{\hat{S}({\gg^L}^*)}^\op$, 
hence we just write $\hat{S}(\gg^*)$ and use the algebra identification,
with the (co)opposite signs $\hat{S}(\gg^*)^{\mathrm{co}}$ or
$\Delta_{\hat{S}(\gg^*)}^\op$ when needed.

As discussed in~\cite{scopr,heisd}, the coproduct is 
equivalently characterized by
\begin{equation}\label{eq:Pbl}
P\blacktriangleright (\hat{f}\hat{g}) = 
m(\Delta_{\hat{S}(\gg^*)}(P)
(\blacktriangleright\otimes\blacktriangleright)(\hat{f}\otimes\hat{g})),
\end{equation}
for all $P\in\hat{S}(\gg^*)$ (for instance, $P = \partial^\mu$) 
and all $\hat{f},\hat{g}\in U(\gg)$. 
Using the action $\blacktriangleright$ we assumed 
that we embedded $\hat{S}(\gg^*)\hookrightarrow
H^R\cong \hat{A}_n$. The right hand version of~(\ref{eq:Pbl}) is
that for all $\hat{u},\hat{v}\in U(\gg^R)$ and $Q\in\hat{S}(\gg^*)$,
\begin{equation}\label{eq:Qbl} 
(\hat{u}\hat{v})\blacktriangleleft Q = 
m((\hat{u}\otimes\hat{v})
(\blacktriangleleft\otimes\blacktriangleleft)\Delta_{\hat{S}(\gg^*)}(Q)).
\end{equation}

\begin{defn} 
The homomorphism 
$\alpha^L: U(\gg^L)\hookrightarrow H^L$ is the inclusion 
$U(\gg^L)\to U(\gg^L)\sharp 1
\hookrightarrow U(\gg^L)\sharp \hat{S}(\gg^*) = H^L$ 
and $\alpha^R:U(\gg^R)\to H^R$ is the inclusion
$\alpha^R:U(\gg^R)\to 1\sharp U(\gg^R)\hookrightarrow 
\hat{S}(\gg^*)\sharp U(\gg^R) = H^R$. Thus, in our writing conventions,
$\alpha^L(\hat{f}) = \hat{f}$ and $\alpha^R(\hat{u}) = \hat{u}$.
Likewise, $\beta^L:U(\gg^L)^\op\to H^L$ and $\beta^R:U(\gg^R)^\op\to H^R$ 
are the unique antihomomorphisms of algebras extending the formulas
(cf.~(\ref{eq:defz}))
\begin{equation}\label{eq:beta}
\begin{array}{l}
\beta^L(\hx_\mu) = \hx_\rho (\OO^{-1})^\rho_\mu = \hy_\mu\, \in H^L.
\\
\beta^R(\hy_\alpha) := \OO^\rho_\alpha \hy_\rho = 
\OO^\rho_\alpha \hx_\sigma (\OO^{-1})^\sigma_\rho = \hat{z}_\alpha \,\in H^R.
\end{array}\end{equation}
\end{defn}
The extension $\beta^L$ exists, because 
the extension of the map $\hx_\mu\mapsto \hy_\mu$ on $\gg$ 
to the antihomomorphism $\beta^L_{T(\gg)}:T(\gg)\to H^L$
maps $[\hx_\alpha,\hx_\beta] - C^\gamma_{\alpha\beta}\hx_\gamma$ 
to $[\hy_\beta,\hy_\alpha]- C^\gamma_{\alpha\beta}\hy_\gamma = 0$;
similarly for $\beta^R$, using~(\ref{eq:zz}). 
\begin{prop}
(i) $H^L$ is a $U(\gg^L)$-bimodule 
via the formula $a.h.b := \alpha^L(a)\beta^L(b)h$,
for all $a,b\in U(\gg^L)$, $h\in H^L$. Likewise, 
$H^R$ is a $U(\gg^R)$-bimodule via $a.h.b := h\beta^R(a)\alpha^R(b)$,
for all $a, b\in U(\gg^R)$, $h\in H^R$.
{\em From now on these bimodule structures are assumed.}

(ii) For any $\hat{f},\hat{g}\in U(\gg^L)$ and any
$\hat{u},\hat{v}\in U(\gg^R)$,

\begin{equation}\label{eq:betabtr}
\beta^L(\hat{g})\btr\hat{f} =\hat{f}\hat{g},\,\,\,\,\,
\hat{u}\btl\beta^R(\hat{v}) =\hat{v}\hat{u}.
\end{equation}
\end{prop}

\begin{pf} (i) The bimodule property of commuting of the left and 
the right $U(\gg^L)$-action is ensured by $[\hx_\mu,\hy_\nu] = 0$.
For the $U(\gg^R)$-actions it boils down to
$[\hat{y}_\mu,\OO^\rho_\nu \hx_\sigma (\OO^{-1})^\sigma_\rho]=0$, 
which follows from the Theorem~\ref{thO}.

(ii) follows from~(\ref{eq:y11x}) and (\ref{eq:z11y}), 
by induction on the filtered degree of $\hat{g}$ 
(respectively, of $\hat{v}$).
\end{pf}
\begin{prop}\label{prop:compHH}
Let $\hat{H}^L := U(\gg^L)\hat\sharp\hat{S}(\gg^*)$ 
and $\hat{H}^R := \hat{S}(\gg^*)\hat\sharp U(\gg^R)$ 
be the completed smash product algebras 
defined in Theorem~\ref{thm:completedsmash}. Then 

(i) the factorwise multiplication $(m\otimes m)(\id\otimes\tau\otimes\id)
: (H^L\otimes H^L)\otimes(H^L\otimes H^L)\to (H^L\otimes H^L)$
(where $\tau$ switches the factors) extends to the unique map
$(H^L\hat\otimes H^L)\otimes(H^L\hat\otimes H^L)\to (H^L\hat\otimes H^L)$
(note that the middle $\otimes$ is not completed!) distributive over formal sums in each of the two $H^L\hat\otimes H^L$-factors. Likewise for $H^R$ in place of $H^L$.

(ii) The inclusions
$H^L\hat\otimes H^L\to\hat{H}^L\hat\otimes\hat{H}^L$, 
$H^R\hat\otimes H^R\to\hat{H}^R\hat\otimes\hat{H}^R$,
$H^L\hat\otimes_{U(\gg^L)} H^L\to\hat{H}^L\hat\otimes_{U(\gg^L)}\hat{H}^L$
and 
$H^R\hat\otimes_{U(\gg^R)} H^R\to\hat{H}^R\hat\otimes_{U(\gg^R)}\hat{H}^R$
are onto; 

(iii) The actions $\btr$ and $\btl$ extend to 
the actions of the completed algebra
$\btr: \hat{H}^L\otimes U(\gg^L)\to U(\gg^L)$ and
$\btl:U(\gg^R)\otimes\hat{H}^R\to U(\gg^R)$.  
\end{prop}
\begin{pf}
(i) The proof is in the vein of 
the proof of Theorem~\ref{thm:completedsmash}.

(ii) The cofiltered components $(H^L)_r = (\hat{H}^L)_r$ agree,
hence both sides of the tensor product inclusions have also
equal cofiltered components. Therefore, the completions are the same.

For (iii) extend the recipe from~(\ref{eq:blact}) and
notice that $\id\hat\sharp\epsilon_S$ kills also all
elements in $U(\gg^L)\hat\sharp\hat{S}(\gg^*)$ not in
$U(\gg^L)\sharp\hat{S}(\gg^*)$ with the result in $U(\gg)$.
On the other hand, {\it there are no completed actions}
$\hat{H}^L\hat\otimes U(\gg^L)\to U(\gg^L)$
and $U(\gg^R)\hat\otimes\hat{H}^R\to U(\gg^R)$
extending $\btr$ and $\btl$.
\end{pf}
\begin{defn} 
The right ideal $I\subset H^L\otimes H^L$ is 
generated by the set of all elements of the form
$\beta^L(\hat{f})\otimes 1 - 1\otimes\alpha^L(\hat{f})$ 
where $\hat{f}\in H^L$.
In other words, $I$ is the kernel of the canonical map 
$H^L\otimes H^L\to H^L\otimes_{U(\gg^L)}H^L$. 

The right ideal $I'\subset H^L\otimes H^L$ 
is the set of all $\sum_i h_i\otimes h'_i \in H^L\otimes H^L$
such that 
$$
\sum_{i,j} (h_i\blacktriangleright \hat{f}_j)
(h'_i\blacktriangleright \hat{g}_j) = 0,
\,\,\,\,\,\mathrm{for}\,\,\,\mathrm{all}\,\, 
\sum_j \hat{f}_j\otimes \hat{g}_j\in U(\gg^L)\otimes U(\gg^L).
$$
Similarly, $\tilde{I} := 
\operatorname{ker}\,(H^R\otimes H^R\to H^R\otimes_{U(\gg^R)}H^R)$ 
is the {\it left} ideal in $H^R\otimes H^R$
generated by all elements of the form
$\alpha^R(\hat{u})\otimes 1 - 1\otimes \beta^R(\hat{u})$, $\hat{u}\in U(\gg^R)$,
and $\tilde{I}'$ is the left ideal in $H^R\otimes H^R$ consisting of all
$\sum_i h_i\otimes h_i'$ such that 
$\sum_{i,j} (\hat{u}_j\btl h_i)(\hat{v}_j\btl h_i') = 0$ 
for all $\sum_j\hat{u}_j\otimes\hat{v}_j\in U(\gg^R)\otimes U(\gg^R)$.
The completions (Appendix A.2) of the ideals 
$I,I'$ and $\tilde{I},\tilde{I}'$ are denoted 
$\hat{I},\hat{I}'\subset H^L\hat\otimes H^L\cong\hat{H}^L\hat\otimes\hat{H}^L$
and 
$\hat{\tilde{I}},\hat{\tilde{I}}'\subset H^R\hat\otimes H^R\cong 
\hat{H}^R\hat\otimes\hat{H}^R$, respectively.
\end{defn}

More generally, for $r\geq 2$, let $I^{(r)}$ be the kernel of the 
canonical projection $(H^L)^{\otimes r}:= H^L
\otimes H^L\otimes \ldots \otimes H^L$ 
($r$ factors) to the tensor product of $U(\gg^L)$-bimodules 
$H^L\otimes_{U(\gg^L)}H^L\otimes_{U(\gg^L)}\ldots \otimes_{U(\gg^L)} H^L$.
$I^{(r)}$ coincides with the smallest right ideal in 
the tensor product algebra $(H^L)^{\otimes r}$ 
which contains $1^{\otimes k}\otimes I \otimes 1^{\otimes (r-k-2)}$
for $k = 0,\ldots,r-2$. Let $I'^{(r)}$ be the set of all elements
$\sum_i h_{1i}\otimes h_{2i}\otimes\ldots h_{ri} \in (H^L)^{\hat\otimes r}$ 
such that for every
$\sum_j u_{1j}\otimes u_{2j}\otimes \ldots 
\otimes u_{rj}\in U(\gg^L)^{\otimes r}$
$$
\sum_{i,j}(h_{1i}\blacktriangleright u_{1j})
        (h_{2i}\blacktriangleright u_{2j})
   \cdots (h_{ri}\blacktriangleright u_{rj}) = 0.
$$

\begin{lem}\label{lemphi}
(i) There is a nondegenerate Hopf pairing 
$$
\langle,\rangle_\phi:U(\gg)\otimes\hat{S}(\gg^*)\to\genfd,
\,\,\,\,\,\,\,\,\langle\hat{u},P\rangle_\phi := \bm\phi_+(\hat{u})(P)(1),
$$
where the action on $1$ is the Fock action 
(on $1$ this amounts to evaluating $\epsilon_{\hat{S}(\gg^*)}$). It
satisfies the \emph{Heisenberg double identity} 
$$
P\btr\hat{u} =
\sum \langle \hat{u}_{(2)}, P\rangle_\phi\hat{u}_{(1)}\,\,\,\,\,
\mbox{for}\,\,\,\,\,P\in\hat{S}(\gg^*)\,\,\,\,\,
\mbox{and}\,\,\,\,\,\hat{u}\in U(\gg).
$$

(ii) For multiindices $J_1,J_2,J,K$ such that $J_1+J_2 = J$, 
$$\bm\phi_+(\hat{x}_K)(\partial^J) = \sum_{K_1+K_2=K} \frac{K!}{K_1!K_2!}
\bm\phi_+(\hat{x}_{K_2})(\partial^{J_1})\bm\phi_+(\hat{x}_{K_1})(\partial^{J_2}).$$

(iii) $\bm\phi_+(\hat{x}_K)(\partial^J)\in\hat{S}(\gg^*)_{|J|-|K|}$ if $|K|<|J|$.

(iv) $\bm\phi_+(\hat{x}_K)(\partial^J)-K!\,
\delta^K_J\in\hat{S}(\gg^*)_1$ if $|K| = |J|$.

(v) For multiindices $K$, $J$ and for the basis 
$\{\partial^K\in S(\gg^*)\}_K$ the identities
$\langle \hat{x}_J,\partial^{K}\rangle_\phi = K!\,\delta^K_J$
hold if $K\geq J$ (in partial order for multiindices),
but in general not otherwise.

(vi) There is a unique family $\{\partial^{\{K\}}\in\hat{S}(\gg^*)\}_K$
which for \underline{all} multiindices $K,J$ satisfies 
$\langle\hat{x}_J,\partial^{\{K\}}\rangle_\phi = K!\, \delta^K_J$.

(vii) Let $f\in\hat{S}(\gg^*)$. Then $\forall r\in\mathbb{N}_0$,
$f_r = \sum_J\frac{1}{J!}\langle
\hat{x}_J,f\rangle_\phi\partial^{\{J\}}_r\in S(\gg)_r$, where the sum 
is finite because $\partial^{\{J\}}_r = 0$ if $r<|J|$. Thus, 
there is a formal sum representation
$f = \underset\longleftarrow\lim{}_r\,f_r 
= \sum_J\frac{1}{J!}\langle\hat{x}_J,f\rangle_\phi\partial^{\{J\}}$.

(viii) $\partial^J = \sum_{|K|\geq |J|} d_{K,J}\partial^{\{K\}}$ 
for some $d_{K,J}\in\genfd$.
\end{lem}
\begin{pf}
(i) is a part of the content of Theorems 3.3 and 3.5 in~\cite{heisd}.

(ii) $\bm\phi_+$ is a right Hopf action, hence the identity follows from the
formula $\Delta(\hat{x}_K) = \sum_{K_1+K_2=K} \frac{K!}{K_1!K_2!}
\hat{x}_{K_1}\otimes\hat{x}_{K_2}$ for the cocommutative coproduct in $U(\gg)$.

(iii) This follows by a simple induction on $|J|-|K|$ using (ii) and 
$\bm\phi_+(1)(\partial^L)=\partial^L\in\hat{S}(\gg^*)$.

(iv) follows by induction on $|K|$ using (ii), (iii) and
$\bm\phi_+(\hat{x}_\mu)(\partial^\mu) = \phi_\mu^\nu$, which by~(\ref{eq:phi}) 
equals $\delta_\mu^\nu$ up to a summand in $\hat{S}(\gg^*)_1$.

(v) This is an application of the formula for $\langle,\rangle_\phi$ in
(i) to the results (iii) and (iv); indeed the elements in $\hat{S}(\gg^*)_1$
vanish when applied to $1$.

(vi) Denote, as in Appendix A.2, by $\pi_r:\hat{S}(\gg^*)\to S(\gg^*)_r$ and
$\pi_{r,r+s}:S(\gg^*)_{r+s}\to S(\gg^*)_r$ the canonical projections.
By~\cite{heisd}, 3.4, the isomorphism $\xi^T:U(\gg)^*\to\hat{S}(\gg^*)$ 
(see~(\ref{eq:xiT}))
of cofiltered algebras identifies the pairing $\langle,\rangle_\phi$
with the evaluation pairing $\langle,\rangle_U:U(\gg)\otimes U(\gg)^*\to\genfd$.
By the properties of $\langle,\rangle_U$, 
for each $r\in\mathbb{N}_0$, the induced 
pairing $\langle,\rangle_r: U(\gg)_r\otimes\,\hat{S}(\gg^*)_r\to\genfd$
characterized by 
$\langle \hat{u},\pi_r(P)\rangle_r = \langle\hat{u},P\rangle_\phi$
for each $\hat{u}\in U(\gg)$, $P\in\hat{S}(\gg^*)$
is {\it nondegenerate}. Thus there is a basis
$\{\partial_r^{\{K\}}\}_{|K|\leq r}$ of the cofiltered component $\hat{S}(\gg^*)_r$ 
dual to the basis $\{\hat{x}_L\}_{|L|\leq r}$ of the filtered component $U(\gg)_r$.
Now 
$\Ker\,\pi_{r,r+s} = \mathrm{Span}\,\{\partial^J,\,r<|J|\leq r+s\}$. 
By (v) $\langle U(\gg)_r,\Ker\,\pi_{r,r+s}\rangle_{r+s} = 0$. 
Therefore for all $K,L$ if $|K|\leq r$, $|L|\leq r$ then  
$\delta^K_L = 
\langle\hat{x}_L,\partial^{\{K\}}_{r+s}\rangle_{r+s}
= \langle \hat{x}_L,\pi_{r,r+s}(\partial^{\{K\}}_{r+s})\rangle_r 
= \langle \hat{x}_L,\partial^{\{K\}}_r\rangle_r$.
By nondegeneracy, 
$\pi_{r,r+s}(\partial^{\{K\}}_{r+s})=\partial^{\{K\}}_r$.
Therefore $\exists!\,\partial^{\{K\}}\in\hat{S}(\gg^*)_{r+s}$ such that
$\pi_r(\partial^{\{K\}}) = \partial^{\{K\}}_r$ for $r\geq |K|$
and $\pi_r(\partial^{\{K\}})=0$ for $r<|K|$; then the requirements
of (vi) hold for $\{\partial^{\{K\}}\}_K$. 

(vii) is now straightforward and (viii) follows from (v) and (vii).
\end{pf}
\begin{theorem}\label{thm:rightIdeal} 
(i) The restriction of 
$\blacktriangleright:H^L\otimes U(\gg^L)\to U(\gg^L)$ 
to $\hat{S}(\gg^*)\otimes U(\gg^L)\to U(\gg^L)$ 
turns $U(\gg^L)$ into a {\bf faithful} 
left $\hat{S}(\gg^*)$-module.

(ii) The right ideals $I,I'$ agree and the left ideals 
$\tilde{I},\tilde{I}'$ agree.

(iii) More generally,
$I^{(r)}=I'^{(r)}$, $\tilde{I}^{(r)}=\tilde{I}'^{(r)}$ for $r\geq 2$.

(iv) Statements (ii) and (iii) hold also for the completed ideals. 
\end{theorem}

\begin{pf} 
We show part (ii) for the right ideals, $I=I'$; 
the method of the proof easily extends to the left ideals, and to
(i), (iii) and (iv).

Let $\sum_\sigma \hat{f}_\sigma\otimes\hat{g}_\sigma\in I$ and 
$v=\hx_{\mu_1}\cdots\hx_{\mu_k}$ a monomial in $U(\gg^L)$. Then 
$$
(\beta^L(v)\btr \hat{f}_\sigma)\hat{g}_\sigma - 
\hat{f}_\sigma \alpha^L(v)\btr\hat{g}_\sigma
= (\hy_{\mu_k}\cdots \hy_{\mu_1}\btr \hat{f}_\sigma) \hat{g}_\sigma
- \hat{f}_\sigma \hx_{\mu_1}\cdots \hx_{\mu_k}\btr\hat{g}_\sigma,
$$
which is zero by Eq.~(\ref{eq:y11x}) and induction on $k$. 
Thus, by linearity, $I\subset I'$.

It remains to show the converse inclusion, $I'\subset I$.
Suppose on the contrary that there is an element 
$\sum_\lambda h_\lambda \otimes h_\lambda'$ in $I'$, but not in $I$; 
then after adding any element in $I$ the sum is still in $I'$ and not in $I$. 
Observe that 
$\hx_J \partial^K\otimes\hx_{J'}\partial^{K'}
= \hx_J \partial^K\otimes\alpha^L(\hx_{J'})\partial^{K'}
= \beta^L(\hx_{J'})\hx_J \partial^K \otimes \partial^{K'}
\,\operatorname{mod}\,I$. 
The tensor factor $\beta(\hx_{J'})\hx_J \partial^K$ belongs to 
$H^L\subset\hat{H}^L$, 
hence it is also a formal linear combination of elements 
of the form $\hx_{J''}\partial^{K''}$.
Therefore, without loss of generality, we can assume
\begin{equation}\label{eq:summ}
\sum_\lambda h_\lambda \otimes h'_\lambda = 
\sum_{J,K,L} a_{JKL} \hx_J\partial^K\otimes \partial^L. 
\end{equation}
Using Lemma~\ref{lemphi} (vi),(vii),(viii) we can 
in~(\ref{eq:summ}) uniquely express
$\partial^K$ as a formal sum in
$\partial^{\lbrace K\rbrace}$ and $\partial^L$ as a formal sum in
$\partial^{\lbrace L\rbrace}$. Therefore,
we can write $\sum_\lambda h_\lambda\otimes h'_\lambda$ as a formal sum  
$$
\sum_\lambda h_\lambda \otimes h'_\lambda = 
\sum_{J,K,L} b_{JKL} \hx_J\partial^{\{K\}}\otimes \partial^{\{L\}},
$$
for some coefficients $b_{JKL}\in\genfd$.
The assumption $\sum_\lambda h_\lambda \otimes h'_\lambda\in I'$ implies
$$
\sum_\lambda (h_\lambda\btr\hx_M)(h'_\lambda\btr\hx_N) = 0.
$$
Choose multiindices $M$ and $N$ such that $(|M|,|N|)$ is a minimal bidegree 
for which $b_{JMN}$ does not vanish for at least some $J$.
By Lemma~\ref{lemphi} (i), the formula 
$\Delta(\hat{x}_M) = \sum_{M_1+M_2=M} \frac{M!}{M_1!M_2!}
\hat{x}_{M_1}\otimes\hat{x}_{M_2}$ for the coproduct in $U(\gg)$,
and Lemma~\ref{lemphi} (vi)
$$
\partial^{\{K\}}\btr\hx_M = 
\sum_{M_1+M_2 = M} {M\choose M_2}
\langle\hx_{M_2},\partial^{\{K\}}\rangle_{\phi}\, \hx_{M_1}
= \left\lbrace\begin{array}{lc}
{M\choose K}\hx_{M-K},&M\geq K\\ 0,&\mathrm{otherwise}.
\end{array}\right.
$$
Therefore, using the minimality of $(|M|,|N|)$,
only the summand with $M = K$ and $N = L$ contributes to the sum and
$$
0 = \sum_\lambda (h_\lambda\btr\hx_M)(h'_\lambda\btr\hx_N) 
= \sum_J b_{JMN}\hat{x}_J,
$$
hence by the linear independence of monomials $\hat{x}_J$, all 
$b_{JMN} = 0$, in contradiction to  
the existence of $J$ with $b_{JMN}$ different from $0$.
\end{pf}

\section{Bialgebroid structures}
\label{sec:bialgebroid}

Let us now use the shorter notation  $\AA^L := U(\gg^L)$,
$\AA^R := U(\gg^R)$. A suggestive symbol $\AA$ 
denotes an abstract algebra in the axioms where either 
$\AA^L$ or $\AA^R$ (or both) may substitute in here intended examples. 
In this section, we equip the isomorphic associative algebras
$H^L$ and $H^R$ with different structures:
$H^L$ is a left $\AA^L$-bialgebroid 
and $H^R$ is a right $\AA^R$-bialgebroid. We start by 
exhibiting the coring structures of these bialgebroids; an $\AA$-coring
is an analogue of a coalgebra where the ground field is replaced
by a noncommutative algebra $\AA$.

\begin{defn}\label{def:coring}
\cite{bohmHbk,BrzWis}
Let $\AA$ be a unital algebra and $C$ an $\AA$-bimodule
with left action $(a,c)\mapsto a.c$ and right action $(c,a)\mapsto c.a$.
A triple $(C,\Delta,\epsilon)$ is an $\AA$-{\bf coring} if 

(i) $\Delta:C\to C\otimes_\AA C$ and $\epsilon:C\to\AA$ 
are $\AA$-bimodule maps; they
are called the coproduct (comultiplication) and the counit;

(ii) $\Delta$ is coassociative: $(\Delta\otimes_{\AA}\id)\circ \Delta 
= (\id\otimes_{\AA}\Delta)\circ\Delta$, where in the codomain 
the associativity isomorphism $(C\otimes_{\AA}C)\otimes_{\AA} C
\cong C\otimes_{\AA}(C\otimes_{\AA} C)$ 
for the $\AA$-bimodule tensor product is understood;

(iii) The counit axioms 
$(\epsilon\otimes_\AA\id)\circ\Delta\cong\id
\cong(\id\otimes_\AA\epsilon)\circ\Delta$
hold, where the identifications of $\AA$-bimodules
$C\otimes_\AA\AA\cong C$, $c\otimes a \mapsto c.a$ and 
$\AA\otimes_\AA C\cong C$, $a\otimes d\mapsto a.d$ are understood.
\end{defn}

\begin{prop}\label{prop:introDelta}
(i) $\exists!$ linear maps 
$\Delta^L:H^L\to H^L\hat\otimes_{\AA^L} H^L$
and $\Delta^R: H^R\to H^R\hat\otimes_{\AA^R} H^R$ such that
$\Delta^L$ and $\Delta^R$ respectively satisfy 
\begin{equation}\label{eq:Pble}
P\blacktriangleright (\hat{f}\hat{g}) = 
m(\Delta^L(P)
(\blacktriangleright\otimes\blacktriangleright)(\hat{f}\otimes\hat{g})),
\,\,\,\,\,\hat{f},\hat{g}\in \AA^L,\,\,P\in H^L,
\end{equation} 
\begin{equation}\label{eq:Qble} 
(\hat{u}\hat{v})\blacktriangleleft Q = 
m((\hat{u}\otimes\hat{v})
(\blacktriangleleft\otimes\blacktriangleleft)\Delta^R(Q)),
\,\,\,\,\,\hat{u},\hat{v}\in \AA^R,\,\,Q\in H^R.
\end{equation}

(ii) $\Delta^L$ is the unique left $\AA^L$-module map 
$H^L\to H^L\hat\otimes H^L$ extending 
$\Delta_{\hat{S}(\gg^*)}:\hat{S}(\gg^*)\to\hat{S}(\gg^*)\hat\otimes\hat{S}(\gg^*)
\subset H^L\hat\otimes H^L$. 
Likewise, $\Delta^R$ is the unique right $\AA^R$-module map  
extending $\Delta_{\hat{S}(\gg^*)}$ from $\hat{S}(\gg^*)$ to $H^R$. 
Equivalently,
\begin{equation}\label{eq:Delta}
\Delta^L (\hat{f}\sharp P)  = \hat{f}\Delta_{\hat{S}(\gg^*)}(P),\,\,\,\,\,\,\,
\Delta^R (Q\sharp\hat{v}) = \Delta_{\hat{S}(\gg^*)}(Q)\hat{v},
\end{equation}
for all $P,Q\in\hat{S}(\gg^*)$, $\hat{f}\in\AA^L$ and $\hat{v}\in\AA^R$.
\newline
In particular, $\Delta^L(\hx_\mu) = \hx_\mu\otimes_{\AA^L}1$
and $\Delta^R(\hy_\mu) = 1\otimes_{\AA^R}\hy_\mu$.

(iii) $\Delta^L(\OO^\mu_\nu) = \OO^\gamma_\nu\otimes_{\AA^L}\OO^\mu_\gamma$,
$\Delta^L (\OO^{-1})^\mu_\nu = 
(\OO^{-1})^\mu_\gamma\otimes_{\AA^L}(\OO^{-1})^\gamma_\nu$,

$\Delta^R (\OO^\mu_\nu) =\OO^\gamma_\nu\otimes_{\AA^R}\OO^\mu_\gamma$,
$\Delta^R (\OO^{-1})^\mu_\nu = 
(\OO^{-1})^\mu_\gamma\otimes_{\AA^R}(\OO^{-1})^\gamma_\nu$,

$\Delta^L(\hy_\nu) = \Delta^L(\hx_\mu(\OO^{-1})^\mu_\nu)
= \hx_\mu (\OO^{-1})^\mu_\gamma\otimes_{\AA^L}(\OO^{-1})^\gamma_\nu
= 1\otimes_{\AA^L}\hy_\nu$,

$\Delta^R(\hx_\nu) = \Delta^R(\hy_\mu\OO^\mu_\nu)
= (1\otimes\hy_\mu)(\OO^\beta_\nu\otimes\OO^\mu_\beta)
= \OO^\beta_\nu\otimes_{\AA^R}\hx_\beta = \hx_\nu\otimes_{\AA^R} 1$.
 
(iv) $(H^L,\Delta^L,\epsilon^L)$ and $(H^R,\Delta^R,\epsilon^R)$
satisfy the axioms for $\AA^L$-coring and $\AA^R$-coring respectively,
provided we replace the tensor product of bimodules 
by the completed tensor of (cofiltered) bimodules and the counit axioms
modify to $(\epsilon\hat\otimes_\AA\id)\circ\Delta\cong j
\cong(\id\hat\otimes_\AA\epsilon)\circ\Delta$ where, 
instead of the identity,
$j$ is the canonical map into the completion (say, $j^L:H^L\hookrightarrow
\hat{H}^L\cong H^L\hat\otimes\genfd\cong\genfd\hat\otimes H^L$).

Taking into account our bimodule structures, 
the counit axioms, Definition~\ref{def:coring} (iii), read 
\begin{equation}\label{eq:epsax}
\begin{array}{l}
\sum\alpha^L(\epsilon^L(h_{(1)}))h_{(2)} = h = 
\sum \beta^L(\epsilon^L(h_{(2)}))h_{(1)},\,\,\,\,\,\,\,h\in H^L\\
\sum h_{(2)}\beta^R(\epsilon^R(h_{(1)})) = h = 
\sum h_{(1)}\alpha^R(\epsilon^R(h_{(2)})),\,\,\,\,\,\,\,h\in H^R.
\end{array}
\end{equation}
(v) The coring structures from (iv) canonically extend to  
an internal $\AA^L$-coring $(\hat{H}^L,\hat\Delta^L,\hat\epsilon^L)$ 
and an internal $\AA^R$-coring $(\hat{H}^R,\hat\Delta^R,\hat\epsilon^R)$  
(see~\cite{bohmInternal}) in the category of complete cofiltered vector spaces
with $\hat\otimes$-tensor product (see Proposition~\ref{prop:compHH} 
and Appendix A.2). Bimodule structures on $\hat{H}^L$, $\hat{H}^R$
involve homomorphisms $\hat\alpha^L:=j^L\circ\alpha^L$,
$\hat\alpha^R:=j^R\circ\alpha^R$, and antihomomorphisms
$\hat\beta^L:=j^L\circ\beta^L$, $\hat\beta^R:=j^R\circ\beta^R$, where
$j^L:H^L\hookrightarrow\hat{H}^L$ and $j^R:H^R\hookrightarrow\hat{H}^R$
are the canonical inclusions. 
\end{prop}

\begin{pf} The equivalence of the two statements in (ii) is evident. 
By Theorem~\ref{thm:rightIdeal} (ii), 
the formulas~(\ref{eq:Pble}) and (\ref{eq:Qble}) 
determine $\Delta^L(P)$ and $\Delta^R(Q)$ uniquely. 
To show the existence, we set the values of $\Delta^L$ and $\Delta^R$ 
by~(\ref{eq:Delta}) and check that~(\ref{eq:Pble}) and~(\ref{eq:Qble}) hold.  
By~(\ref{eq:Pbl}) and~(\ref{eq:Qbl}) we already know this for
$P,Q\in\hat{S}(\gg^*)$. Using the
action axiom for $\blacktriangleright$, observe that
$$
\begin{array}{lcl}
\hx_\mu\blacktriangleright (P
\blacktriangleright (\hat{f}\hat{g})) &=& 
\hx_\mu \cdot m(\Delta^L(P)
(\blacktriangleright\otimes\blacktriangleright)(\hat{f}\otimes\hat{g}))
\\&=& m(\hx_\mu \Delta^L(P)
(\blacktriangleright\otimes\blacktriangleright)(\hat{f}\otimes\hat{g}))
\\&\overset{(\ref{eq:Delta})}=& m(\Delta^L(\hx_\mu P)
(\blacktriangleright\otimes\blacktriangleright)(\hat{f}\otimes\hat{g}))
\end{array}
$$
for all $\hat{f},\hat{g}\in\AA^L$, 
hence~(\ref{eq:Pble}) holds for all $P\in H^L$. 
Likewise check~(\ref{eq:Qble}) for all $Q\in H^R$.
Conclude (i). The statement in
(ii) that $\Delta^L$, $\Delta^R$ extend $\Delta_{\hat{S}(\gg^*)}$,
 $\Delta_{\hat{S}(\gg^*)}^\op$ is the statement that
(\ref{eq:Pble}),(\ref{eq:Qble}) specialize to 
(\ref{eq:Pbl}),(\ref{eq:Qbl}) when $P,Q\in\hat{S}(\gg^*)$. The rest of
(ii) follows from uniqueness in (i).

(iii) By Theorem~\ref{thm:rightIdeal} (ii), 
the first 4 formulas follow from~(\ref{eq:DeltaO}),(\ref{eq:DeltamO}),
(\ref{eq:DeltaRO}),(\ref{eq:DeltaRmO}). 
The formulas for $\Delta^L(\hy_\alpha)$
and $\Delta^R(\hx_\alpha)$ are straightforward.

(iv) To show that $\Delta^L$ is an $\AA^L$-bimodule map
note that by (ii) $\Delta^L$ commutes with the left $\AA^L$-action.
It remains to show that $\Delta^L$ also commutes  
with the right $\AA^L$-action. This is sufficient to check on the
generators $\hx_\mu$ of $\AA^L$ and arbitrary $P\in\hat{S}(\gg^*)$:

$$\begin{array}{lcl}
\Delta^L(P).\hx_\mu &=& \sum P_{(1)}\otimes_{\AA^L}\beta(\hx_\mu) P_{(2)}
\\ &=& \sum P_{(1)}\otimes_{\AA^L}\alpha(\hx_\nu)(\OO^{-1})^\nu_\mu P_{(2)}
\\ &=& \sum \beta(\hx_\nu) P_{(1)}\otimes_{\AA^L}(\OO^{-1})^\nu_\mu P_{(2)}
\\ &=& \sum \hx_\gamma(\OO^{-1})^\gamma_\nu P_{(1)}\otimes_{\AA^L}(\OO^{-1})^\nu_\mu P_{(2)}
\\ &=& \Delta^L(\hat{x}_\nu(\OO^{-1})^\nu_\mu P)
\\ &=& \Delta^L(\beta(\hat{x}_\mu)(P))
\end{array}$$

By Theorem~\ref{thm:rightIdeal} (iii) for $r=3$, 
the action axiom for $\blacktriangleright$ and associativity in $H^L$
implies the coassociativity of $\Delta^L$. 

We exhibit the counits $\epsilon^L$ and $\epsilon^R$ 
(and their completed versions 
$\hat\epsilon^L:\hat{H}^L\to\AA^L$, 
$\hat\epsilon^R:\hat{H}^R\to\AA^R$) by the corresponding actions on $1$,
\begin{equation}\label{eq:epsdef}
\epsilon^L(h) := h\blacktriangleright 1_{\AA^L},\,\,\,\,\,\,\,
\epsilon^R(h) := 1_{\AA^R}\blacktriangleleft h.
\end{equation}
The counit axioms~(\ref{eq:epsax}) for $\epsilon^L$ 
are checked on the generators $\hx_\mu$:
\begin{eqnarray*}
\sum \alpha(\epsilon^L(\hx_{\mu(1)}))\hx_{\mu(2)}
= \alpha(\epsilon^L(\hx_\mu))1 = \hx_\mu,
\\
\sum \beta(\epsilon^L(\hx_{\mu(2)}))\hx_{\mu(1)}
= \beta(\epsilon^L(1))\hx_\mu = \hx_\mu. 
\end{eqnarray*} 
\nodo{
or, if we use the form $\Delta(\hx_\mu)=\mathcal{O}^\nu_\mu\otimes\hx_\nu$,
\begin{eqnarray*}
\sum \beta(\epsilon^L(\hx_{\mu(2)}))\hx_{\mu(1)} 
= \beta(\epsilon^L(\hx_\nu))\mathcal{O}^\nu_\mu 
=\hx_\lambda(\mathcal{O}^{-1})^\lambda_\nu\mathcal{O}^\nu_\mu = \hx_\mu,
\\
\sum \alpha(\epsilon^L(\hx_{\mu(1)}))\hx_{\mu(2)}
=\alpha(\epsilon^L(\mathcal{O}^\nu_\mu))\hx_\nu 
= \alpha(\delta^\nu_\mu 1)\hx_\nu  = \hx_\mu.
\end{eqnarray*}
}
Similarly, one checks the counit identities for $\epsilon^R$. 

Using formal expressions in the completions, (v) is straightforward.
\end{pf}

\begin{defn} 
(Modification of~\cite{bohmHbk,Bohm,BrzMilitaru}, cf.~\ref{def:leftbialg}).
Given an algebra $\AA$, a {\bf formally completed left $\AA$-bialgebroid} 
$(H,m,\alpha,\beta,\Delta,\epsilon)$ consists of the following data.
$H$ is a cofiltered vector space and $(H,m)$ an associative algebra
with multiplication $m$ distributive with respect to 
the formal sums in each argument (Appendix A.2) and the factorwise
multiplication on $H\otimes H$ extends to a multiplication
$(H\hat\otimes H)\otimes (H\hat\otimes H)$ distributive with respect to
the formal sums in each argument;
$\alpha:\AA\to H$ and $\beta: \AA^\op\to H$ 
are fixed algebra homomorphisms with commuting images;
$H$ is equipped with a structure of an $\AA$-bimodule via the formula
$a.h.a':= \alpha(a)\beta(a')h$; 
$\Delta : H\to H\hat\otimes_{\AA}H$ is an $\AA$-bimodule map, 
coassociative and with counit $\epsilon:H\to\AA$ understood with
respect to the completed tensor product $\hat\otimes$
and the counit axiom modifies to  
$(\epsilon\hat\otimes_\AA\id)\circ\Delta\cong j
\cong(\id\hat\otimes_\AA\epsilon)\circ\Delta$ where 
$j:H\to H\hat\otimes\genfd\cong\hat{H}\cong\genfd\hat\otimes H$ 
is the canonical map into the completion;
both $\Delta$ and $\epsilon$ should be distributive 
with respect to formal sums. It is required that

(i) $\epsilon$ is a {\bf left character} on the $\AA$-ring
$(H,m,\alpha)$ in the sense that the formula 
$h\otimes\hat{f}\mapsto \epsilon(h\alpha(\hat{f}))$ defines an action 
$H\otimes\AA\to\AA$ extending the left regular action 
$\AA\otimes\AA\to\AA$;

(ii) the coproduct $\Delta :H\to H\hat\otimes_{\AA}H$ 
corestricts to the {\bf formal Takeuchi product}
$$
H\hat\times_{\AA}H\subset H\hat\otimes_{\AA}H
$$
which is by definition the $\AA$-subbimodule, consisting of
all formal sums $b = \sum_\lambda b_\lambda\otimes_\AA b'_\lambda$,
where also 
$\sum_\lambda b_\lambda\otimes_\genfd b'_\lambda\in H\hat\otimes_\genfd H$ 
is formal and such that
$$
\sum_\lambda b_\lambda\otimes_\AA b'_\lambda \alpha(a) = \sum_\lambda b_\lambda \beta(a)\otimes_\AA b'_\lambda, \,\,\forall a\in \AA.
$$
(iii) The corestriction $\Delta| :H\to H\hat\times_{\AA}H$ 
is an algebra map.
\end{defn}

Notice that $H\hat\otimes_\AA H$
does not carry a well-defined multiplication induced from $H\hat\otimes H$,
unlike $H\hat\times_\AA H$ which does. This explains the need for (ii). 
Indeed, (ii) implies that for any $b\in H\hat\times H$ and for any formal sum 
$c = \sum_\mu c_\mu\otimes_\AA c'_\mu\in H\hat\otimes_\AA H$ 
the product $b\cdot c$ obtained by lifting $b$ and $c$ 
to $H\hat\otimes_\genfd H$ is 
a well defined element of $H\hat\otimes_\AA H$ 
(does not depend on the lifting as a formal sum in $H\hat\otimes_\genfd H$); 
if moreover $c\in H\hat\times_\AA H$ then $b\cdot c\in H\hat\times_\AA H$. 
Thus $H\hat\times_\AA H$ is an algebra and (iii) makes sense.

Interchanging the left and right sides 
in all modules and binary tensor products 
in the definition of a left $\AA$-bialgebroid, we get
a {\bf right $\AA$-bialgebroid}~(\cite{bohmHbk}). The $\AA$-bimodule
structure on $H$ is then given by $a.h.b := h\alpha(b)\beta(a)$. 
In short, $(H,m,\alpha,\beta,\Delta,\epsilon)$ is
a right $\AA$-bialgebroid iff $(H,m,\beta,\alpha,\Delta^\op,\epsilon)$ 
is a left $\AA^\op$-bialgebroid; analogously with the completed versions.  

Let us return to our candidate examples $H^L$ and $H^R$. Regarding
that $H^L\hat\otimes H^L=\hat{H}^L\hat\otimes\hat{H}^L$ it may be
convenient to have all modules completed to start with, hence
considering the completed smash product algebras 
$\hat{H}^L$ and $\hat{H}^R$ (Theorem~\ref{thm:completedsmash} 
and Proposition~\ref{prop:compHH}). One of the advantages is
that for $\hat{H}^L$ and $\hat{H}^R$ 
the internal coring axiom is not modified for $\hat\Delta^L$.
Still, the expectation that
all objects and morphisms are in the completed sense is {\it not} true, 
as the multiplication $m:H^L\otimes H^L\to H^L$ 
can be extended to $\hat{H}^L\otimes\hat{H}^L\to\hat{H}^L$ 
but can not be extended to a function on the completed tensor product
$\hat{H}^L\hat\otimes\hat{H}^L\to\hat{H}^L$  
distributive over formal sums.
The thesis~\cite{StojicPhD} alternatively introduces 
a canonical tensor product
on a more complicated category of filtrations of cofiltrations;
it involves less drastic completions in general and 
admits a truly internal bialgebroid structure on $H^L$.

\begin{prop} \label{prop:Hisbialg}
$(H^L,m,\alpha^L,\beta^L,\Delta^L,\epsilon^L)$ and
$(\hat{H}^L,\hat{m},\hat\alpha^L,\hat\beta^L,\hat\Delta^L,\hat\epsilon^L)$ 
have a structure of formally completed left $\AA^L$-bialgebroids.
Likewise, 
$(H^R,m,\alpha^R,\beta^R,\Delta^R,\epsilon^R)$ and
$(\hat{H}^R,\hat{m},\hat\alpha^R,\hat\beta^R,\hat\Delta^R,\hat\epsilon^R)$
are formally completed right $\AA^R$-bialgebroids. 
\end{prop}

\begin{pf} 
The coring axioms are checked in Proposition~\ref{prop:introDelta}. 

To check that the rule $\sum_\lambda h_\lambda\otimes\hat{f}_\lambda\mapsto
\sum_\lambda\epsilon^L(h_\lambda\alpha(\hat{f}_\lambda))$ (for finite sums)
is an action and (i) holds for $\epsilon^L$, observe from 
the definition~(\ref{eq:epsdef}) that
$\epsilon^L(h\alpha(\hat{f})) =  
h\alpha(\hat{f})\blacktriangleright 1 = 
h\blacktriangleright\hat{f}$, for all $\hat{f}\in\AA^L$, $h\in H^L$.
Analogously check (i) for $\epsilon^R,\hat\epsilon^L,\hat\epsilon^R$. 

To show that $\hat\Delta^L$ corestricts 
to the formal Takeuchi product $\hat{H}^L\hat\times_{\AA^L}\hat{H}^L$, 
calculate for $P\in\hat{H}^L$ and $\hat{f},\hat{g},\hat{h}\in\AA^L$,
$$\begin{array}{lcl}
((P_{(1)}\hat\beta^L(\hat{g})\btr\hat{f})\cdot (P_{(2)}\btr\hat{h})&=&
(P_{(1)}\btr(\hat\beta^L(\hat{g})\btr\hat{f}))\cdot
(P_{(2)}\btr\hat{h}) 
\\ &\overset{(\ref{eq:betabtr})}=& (P_{(1)}\btr (\hat{f}\hat{g}))
\cdot (P_{(2)}\btr\hat{h})
\\ &=& P\btr (\hat{f}\hat{g}\hat{h})
\\ &=& (P_{(1)}\btr\hat{f})\cdot((P_{(2)}\hat\alpha^L(\hat{g}))\btr\hat{h}),
\end{array}$$
thus, by Theorem~\ref{thm:rightIdeal} (ii), (iv),  
$P_{(1)}\hat\beta^L(\hat{g})\otimes_{\AA^L} P_{(2)}
= P_{(1)}\otimes_{\AA^L} P_{(2)}\hat\alpha^L(\hat{g})$,
hence $\hat\Delta^L(P)\in\hat{H}^L\hat\times_{\AA^L}\hat{H}^L$.  
 
We now check directly that the corestriction 
$\hat\Delta^L:\hat{H}^L\to \hat{H}^L\hat\times_{\AA^L}\hat{H}^L$ is
a homomorphism of algebras, 
$$
\hat\Delta^L(h_1 h_2) = \hat\Delta^L(h_1)\hat\Delta^L(h_2)
\,\,\,\,\,\,\,\mbox{for all }
h_1,h_2\in H^L.
$$
To this aim, recall that $\Delta_{\hat{S}(\gg^*)}:\hat{S}(\gg^*)
\to\hat{S}(\gg^*)\hat\otimes \hat{S}(\gg^*)$ is a homomorphism,
and that by Proposition~\ref{prop:introDelta} (ii),
$\hat{\Delta}^L|_{1\sharp \hat{S}(\gg^*)}$ is the composition
$$
1\sharp \hat{S}(\gg^*)\cong \hat{S}(\gg^*)
\overset{\Delta_{\hat{S}(\gg^*)}}\longrightarrow 
\hat{S}(\gg^*)\hat\otimes \hat{S}(\gg^*)\hookrightarrow
\hat{H}^L\hat\times_{\AA^L}\hat{H}^L,$$
hence homomorphism as well (the inclusion is a homomorphism, 
because the product is factorwise).
We use this when applying to the tensor factor 
$P_{(2)} Q$ in the calculation
$$\begin{array}{lcl}
\hat\Delta^L((u\sharp P)(v\sharp Q)) &=& 
\hat\Delta^L(u (P_{(1)}\blacktriangleright v) \sharp P_{(2)} Q) 
\\
&=& (u(P_{(1)}\blacktriangleright v)\sharp P_{(2)} Q_{(1)}) 
\otimes (1\sharp P_{(3)}Q_{(2)}).
\\
&=& [(u\sharp P_{(1)})(v\sharp Q_{(1)})]\otimes (1\sharp P_{(2)} Q_{(2)})
\\
&=& (u (P_{(1)}\blacktriangleright v)\sharp P_{(2)} Q_{(1)})
\otimes(1\sharp P_{(3)} Q_{(2)}).
\\
&=& [(u\sharp P_{(1)})\otimes(1\sharp P_{(2)})][(v\sharp Q_{(1)})\otimes(1\sharp Q_{(2)})] 
\\ &=& 
\hat\Delta^L(u\sharp P)\hat\Delta^L(v\sharp Q).
\end{array}
$$\end{pf}

\section{The antipode and Hopf algebroid}
\label{sec:antipode}

A Hopf algebroid is roughly a bialgebroid with an antipode. 
In the literature, there are several nonequivalent versions.
In the framework of G. B\"ohm~\cite{bohmHbk}, 
there are two variants which are equivalent 
if the antipode is bijective (as it is here the case): 
nonsymmetric and symmetric. The
{\it nonsymmetric} involves one-sided bialgebroid with an antipode map
satisfying axioms which involve both the antipode map and its
inverse. The {\it symmetric} version involves two bialgebroids
and axioms neither involve nor require the inverse of the antipode. 
We choose this version here, 
because we naturally constructed two actions, 
$\blacktriangleright$ and $\blacktriangleleft$,
which lead to the two coproducts, $\Delta^L$ and $\Delta^R$, 
as exhibited in Section~\ref{sec:bialgebroid}. 

\begin{defn}
Given two algebras $\AA^L$ and $\AA^R$ with
fixed isomorphism $(\AA^L)^{\mathrm{op}}\cong\AA^R$, 
a {\bf symmetric Hopf algebroid}~(\cite{bohmHbk}) is a pair of 
a left $\AA^L$-bialgebroid $H^L$ and a
right $\AA^R$-bialgebroid $H^R$, isomorphic
and identified as algebras $H\cong H^L\cong H^R$, 
such that the compatibilities
\begin{equation}\label{eq:alphaepsilonbeta}
\begin{array}{lr}
\alpha^L\circ\epsilon^L\circ\beta^R = \beta^R,
&
\beta^L\circ\epsilon^L\circ\alpha^R = \alpha^R,
\\
\alpha^R\circ\epsilon^R\circ\beta^L = \beta^L,
&
\beta^R\circ\epsilon^R\circ\alpha^L = \alpha^L,
\end{array}\end{equation}
hold between the source and target maps 
$\alpha^L,\alpha^R,\beta^L,\beta^R$, 
and the counits $\epsilon^L,\epsilon^R$; 
the comultiplications $\Delta^L$ and $\Delta^R$ 
satisfy the compatibility relations
\begin{equation}\label{eq:deltalr1}
(\Delta^R\otimes_{\AA^L}\id)\circ\Delta^L = 
(\id\otimes_{\AA^R}\Delta^L)\circ\Delta^R
\end{equation}
\begin{equation}\label{eq:deltalr2}
(\Delta^L\otimes_{\AA^R}\id)\circ\Delta^R = 
(\id\otimes_{\AA^L}\Delta^R)\circ\Delta^L
\end{equation}
and there is a map $\mcS:H\to H$, called the {\bf antipode} 
which is an antihomomorphism of algebras and satisfies 
\begin{equation}\label{eq:antipodeMain}
\begin{array}{c}
\mcS\circ \beta^L =\alpha^L,\,\,\,\,\,\,\,\,\,\,\, 
\mcS\circ\beta^R = \alpha^R
\\
m\circ(\mcS\otimes\id)\circ\Delta^L = \alpha^R\circ\epsilon^R
\\
m\circ(\id\otimes \mcS)\circ\Delta^R = \alpha^L\circ\epsilon^L
\end{array}
\end{equation}
A {\bf formally completed symmetric Hopf algebroid} is defined
analogously as a pair of left and right formally completed bialgebroid
with antipode $\mcS$ satisfying~(\ref{eq:antipodeMain})
and the compatibilities~(\ref{eq:alphaepsilonbeta}),(\ref{eq:deltalr1}),(\ref{eq:deltalr2}) 
satisfied with the tensor products
replaced with the completed ones.
\end{defn}

\begin{theorem}
Data $\AA^L = U(\gg^L)$, $\AA^R = U(\gg^R)$ together with either

(i) $\hat{H}^L:=U(\gg^L)\hat\sharp \hat{S}(\gg^*)$, 
$\hat{H}^R:=\hat{S}(\gg^*)\hat\sharp U(\gg^R)$, 
$\hat\epsilon^L,\hat\epsilon^R,\hat\alpha^L,\hat\beta^L,
\hat\alpha^R,\hat\beta^R$
from Section~\ref{sec:bimodule} and $\hat\Delta^L$, $\hat\Delta^R$ 
defined in Section~\ref{sec:bialgebroid},

(ii) or $H^L = U(\gg^L)\sharp \hat{S}(\gg^*)$, 
$H^R:=\hat{S}(\gg^*)\sharp U(\gg^R)$, 
$\epsilon^L,\epsilon^R,\alpha^L,\beta^L,\alpha^R,\beta^R$
from Section~\ref{sec:bimodule} and $\Delta^L$, $\Delta^R$ 
defined in Section~\ref{sec:bialgebroid},

\noindent form a formally completed symmetric Hopf algebroid. 
The antipode map for (i) is the unique 
homomorphism of algebras $\mcS:\hat{H}\to\hat{H}$ 
distributive over formal sums and such that
$$\mcS(\partial^\nu) = -\partial^\nu,$$
and the antipode for (ii) is 
its restriction $\mcS = \mcS|:H\to H$.
The antipode $\mcS$ is bijective in both cases, 
and by distributivity over formal sums it follows that
$\mcS(\OO) = \mcS (e^\CC) = e^{-\CC} = \OO^{-1}$ and
\begin{equation}\label{eq:sshy}
\mcS(\hy_\mu)=\hx_\mu.
\end{equation}
For a general Lie algebra $\gg$, $\mcS^2\neq\id$. 
More precisely,
\begin{equation}\label{eq:sshx}
\mcS^2(\hy_\mu)=\mcS(\hx_\mu) = \hy_\mu - C^\lambda_{\mu\lambda},\,\,\,\,\,
\mcS^{-2}(\hx_\mu) = \mcS^{-1}(\hy_\mu) = \hx_\mu - C^\lambda_{\mu\lambda}
\end{equation}
\begin{equation}\label{eq:ss2}
\mcS^2(\hx_\mu) = \hx_\mu + C^\lambda_{\mu\lambda},\,\,\,\,\,\,\,
\mcS^{-2}(\hy_\mu) = \hy_\mu + C^\lambda_{\mu\lambda}.
\end{equation}
with the summation over $\lambda$ understood.
\end{theorem}
\begin{pf}
In this proof, we simply write $\epsilon^L,\Delta^L$ etc. without
hat symbol, as it is not essential for the arguments below
which work for both versions. We proved that the above
data give bialgebroids (Proposition~\ref{prop:Hisbialg}).  

One checks the relations~(\ref{eq:alphaepsilonbeta}) 
on generators, for which 
$\alpha^R(\hy_\mu)=\hy_\mu$, $\beta^R(\hy_\mu)=\OO^\rho_\mu \hy_\rho 
= \OO^\rho_\mu\hy_\sigma(\OO^{-1})^\sigma_\rho$,
$\alpha^L(\hx_\mu)=\hx_\mu$, $\beta^L(\hx_\mu)=\hy_\mu$.

Regarding that $\Delta^L$ and $\Delta^R$ restricted to 
$\hat{S}(\gg^*)$ coincide with $\Delta_{\hat{S}(\gg^*)}$, 
(\ref{eq:deltalr1}) and (\ref{eq:deltalr2}) restricted to $\hat{S}(\gg^*)$ 
reduce to the coassociativity. 
Algebra $H^L$ is generated by $\hat{S}(\gg^*)$ and $\gg^L$, 
so it is enough to check~(\ref{eq:deltalr1}),(\ref{eq:deltalr2}) 
also on $\hy_\mu=\hx_\nu(\OO^{-1})^\nu_\mu$. 
This follows from the matrix identities
$$
\Delta^L(\hx\OO^{-1}) = \hx\OO^{-1}\otimes_{\AA^L}\OO^{-1}
= \OO^{-1}\OO\otimes_{\AA^L}\hx\OO^{-1} = 1\otimes_{\AA^L}\hy,
$$
$$
\Delta^R(\hy) = 1\otimes_{\AA^R}\hy = \hy\otimes_{\AA^R}\OO^{-1}
= \hx\OO^{-1}\otimes_{\AA^R}\OO^{-1}.
$$

Formula $\mcS(\partial^\mu) = -\partial^\mu$ clearly extends to a unique
continuous antihomomorphism of algebras on the formal power series
ring $\hat{S}(\gg^*)$. Similarly, by functoriality of
$\gg\mapsto U(\gg)$, the antihomomorphism of Lie algebras,
$\mcS:\gg^R\to\gg^L$, $\hy_\mu\mapsto \hx_\mu$, 
extends to a unique antihomomorphism $U(\gg^R)\to U(\gg^L)$.  
Regarding that $U(\gg^R)$ and $\hat{S}(\gg^*)$ generate $H^R$,
it is sufficient to check that $\mcS$ is compatible
with the additional relations in the smash product, namely
$[\partial^\mu,\hy_\nu] = \left(\frac{\CC}{e^{\CC}-1}\right)^\mu_\nu$. 
Then $\mcS([\partial^\mu,\hx_\nu]) = 
\mcS\left(\frac{-\CC}{e^{-\CC}-1}\right)^\mu_\nu
= \left(\frac{\CC}{e^\CC-1}\right)^\mu_\nu = 
\left(e^{-\CC} \frac{-\CC}{e^{-\CC}-1}\right)^\mu_\nu$, which
equals $(e^{-\CC})^\rho_\nu [\hx_\rho, -\partial^\mu]
= [\mcS(\hy_\rho\OO^\rho_\nu), -\partial^\mu]
= [\mcS(\hx_\nu),\mcS(\partial^\mu)]$.

To exhibit the inverse $\mcS^{-1}$, we similarly check that the obvious
formulas $\mcS^{-1}(\hx_\mu) = \hy_\mu$, $\mcS^{-1}(\partial^\mu)=\partial^\mu$
define a unique continuous (that is, distributive over formal sums)
antihomomorphism $\mcS^{-1}:H\to H$.

For~(\ref{eq:sshx}) calculate
$\mcS(\hx_\mu) = \mcS(\hy_\rho\OO^\rho_\mu)
= \mcS(\OO^\rho_\mu)\mcS(\hy_\rho) = (\OO^{-1})^\rho_\mu \hx_\rho
= (\OO^{-1})^\rho_\mu \hy_\sigma \OO^\sigma_\rho$ and use 
$[\OO^\rho_\mu,\hy_\sigma]= -C_{\tau\sigma}^\rho\OO^\tau_\mu$
in the last step. Similarly, we get $\mcS^{-1}(\hy_\mu) = 
\OO^\rho_\mu \hx_\sigma (\OO^{-1})^\sigma_\rho$
and use $[\OO^\rho_\mu,\hx_\sigma]= -C_{\tau\sigma}^\rho\OO^\tau_\mu$ for 
the second formula in~(\ref{eq:sshx}). 
Notice that $\mcS^{-1}(\hy_\mu) = \hat{z}_\mu$ from
Theorem~\ref{th:btl}, formula~(\ref{eq:z11y}).
For~(\ref{eq:ss2}) similarly use the matrix identities 
$\mcS^2(\hx) = \mcS(\OO^{-1}\hy\OO)= \OO^{-1} \hx \OO$,
$\mcS^{-2}(\hy) = \OO \hy\OO^{-1}$. 

The formula $\mcS(\beta^L(\hx_\mu))=\mcS(\hy_\mu)=\hx_\mu=\alpha^L(\hx_\mu)$ shows 
$\mcS\circ \beta^L =\alpha^L$ for the generators of $\AA^L$.
Likewise for the rest of the identities~(\ref{eq:antipodeMain}).
\end{pf}

\section{Conclusion and perspectives.}

We have equipped the noncommutative phase spaces of Lie algebra type 
with the structure of a version of a Hopf algebroid over $U(\gg)$.
That roughly means that we have found a left $U(\gg)$-bialgebroid $H^L$, 
and a right $U(\gg)^{\mathrm{op}}$-bialgebroid $H^R$, 
which are canonically
isomorphic as associative algebras $H^L\cong H^R$, and
an antipode map $\mcS$ satisfying a number of axioms involving 
a completed tensor product $\hat\otimes$.

Hopf algebroids allow a version of Drinfeld's twisting cocycles 
studied earlier in the context of deformation quantization~(\cite{xu}),
and are a promising tool for extending many constructions
to the noncommutative case, and a planned
direction for our future work. One can find a cocycle 
which can be used to twist 
the Hopf algebroid corresponding to the abelian Lie algebra
(i.e. the Hopf algebroid structure on 
the completion of the usual Weyl algebra)
to recover the Hopf algebroid of the phase space for any
other Lie algebra of the same dimension~(\cite{twist}). 
More importantly for applications,
along with the phase space one can systematically twist 
many geometric structures, including differential forms, 
from the undeformed to the deformed case.
This has earlier been studied in
the case of $\kappa$-spaces (e.g. in~\cite{tajronkov}), 
while the work for general finite-dimensional Lie algebras (and 
for some nonlinear star products) is in progress. 

\vskip .1in {\footnotesize {\bf Acknowledgments.} 
We thank L. El Kaoutit, V. Roubtsov, T.~Brzezi\'nski, 
T.~Maszczyk and G. B\"ohm
for reading fragments of this work and advice. We thank 
A.~Borowiec, J. Lukierski and A. Pacho\l\ for discussions. 
S.M. has been supported by Croatian Science Foundation, 
project IP-2014-09-9582.
A part of the work has been done at IRB, Zagreb, 
and a substantial progress has also been made during 
the visit of Z.\v{S}. to IH\'ES in November 2012. 
Z.\v{S}. thanks G. Garkusha for the invitation to present a related 
{\it Algebra and topology seminar} at Swansea, July 15, 2014.
}

\appendix

\subsection*{A.1 Commutation $[\hx_\alpha,\hy_\beta] = 0$}

\begin{prop} \label{prop:xycommute}
The identity $[\hx_\mu,\hy_\nu] = 0$ holds 
in the realization $\hx_\mu = x_\sigma \phi^\sigma_\mu=
x_\sigma\left(\frac{-\CC}{e^{-\CC}-1}\right)^\sigma_\mu$,
$\hy_\mu = x_\rho \tilde\phi^\rho_\mu =
x_\rho\left(\frac{\CC}{e^{\CC}-1}\right)^\rho_\mu$, where
$\CC^\mu_\nu = C^\mu_{\nu\gamma}\partial^\gamma$ 
(cf. the equations ~(\ref{eq:parhx},\ref{CCmatrix},\ref{eq:phi})).
\end{prop}
\begin{pf}
For any formal series $P = P(\partial)$ in $\partial$-s, 
$[P,\hx_\mu] = \frac{\partial P}{\partial (\partial^\mu)} =:\delta_\mu P$. 
In particular (cf.~\cite{ldWeyl}), 
from $[\hx_\mu,\hx_\nu] = C_{\mu\nu}^\lambda\hx_\lambda$,
one obtains a formal differential equation
for $\phi^\sigma_\mu$,
\begin{equation}\label{eq:formalphi}
(\delta_\rho \phi^\gamma_\mu)\phi^\rho_\nu - (\delta_\rho
\phi^\gamma_\nu)\phi^\rho_\mu = C^\sigma_{\mu\nu}
\phi^\gamma_\sigma.
\end{equation}
By symmetry $C^i_{jk}\mapsto -C^i_{jk}$ 
the same equation holds 
with $(-\tilde\phi) = \frac{-\CC}{e^{\CC}-1}$ in the place of $\phi$. 
Similarly, the equation $[\hx_\mu,\hy_\nu] = 0$, i.e.
$[x_\gamma \phi^\gamma_\mu,x_\beta \tilde\phi^\beta_\nu] = 0$, is
equivalent to
\begin{equation}\label{eq:formalphipsi}
(\delta_\rho\phi^\gamma_\mu)\tilde\phi^\rho_\nu -
(\delta_\rho\tilde\phi^\gamma_\nu)\phi^\rho_\mu = 0
\end{equation}
Recall that $\phi = \frac{-\CC}{e^{-\CC}-1} = \sum_{N = 0}^\infty
(-1)^N \frac{B_N}{N!} (\CC^N)^i_j$, where $B_N$ are the Bernoulli
numbers, which are zero unless $N$ is either even or $N=1$. 
Hence
$\tilde\phi = \frac{\CC}{e^{\CC}-1} = 
\sum_{N = 0}^\infty \frac{B_N}{N!}\CC^N =
\frac{B_1}{2}\CC+\sum_{N\,\mbox{even}}^\infty \frac{B_N}{N!}\CC^N$ 
and $\phi-\tilde\phi = -2\frac{B_1}{2}\CC = \CC$.
Notice that $\frac{\partial \CC^\alpha_\beta}{\partial
(\partial^\mu)} = C^\alpha_{\beta\mu}$. Therefore, subtracting
(\ref{eq:formalphipsi}) from (\ref{eq:formalphi}) gives
the condition
$$
(\delta_\rho \phi^\gamma_\mu)\CC^\rho_\nu - C^\gamma_{\nu\rho}
\phi^\rho_\mu = C^\sigma_{\mu\nu} \phi^\gamma_\sigma.
$$
$\CC$ is homogeneous of degree $1$ in $\partial^\mu$-s, 
so we can split this condition into the parts 
of homogeneity degree $N$:
\begin{equation}\label{eq:CCN}
[\delta_\rho (\CC^N)^\gamma_\mu]\CC^\rho_\nu -
(\delta_\rho\CC^\gamma_\nu)(\CC^N)^\rho_\mu =
C^\sigma_{\mu\nu}(\CC^N)^\gamma_\sigma,
\end{equation}
where the overall factor of $(-1)^N B_N/N!$ has been taken out.
Hence the proof is reduced to the following lemma:\end{pf}
\begin{lem}
The identities~(\ref{eq:CCN}) hold for $N = 0,1,2,\ldots$.
\end{lem}
\begin{pf} For $N = 0$, (\ref{eq:CCN}) 
reads $C^\gamma_{\nu\mu} = C^\gamma_{\mu\nu}$, 
which is the antisymmetry of the bracket. For $N=1$ it follows
from the Jacobi identity:
$$
(C^\gamma_{\mu\rho}C^\rho_{\nu\tau}-C^\gamma_{\nu\rho}C^\rho_{\mu\tau})\partial^\tau
= C^\rho_{\mu\nu}C^\gamma_{\rho\tau}\partial^\tau.
$$
Suppose now (\ref{eq:CCN}) holds for given $N=K\geq 1$. Then
$$
C^\gamma_{\mu\nu}(\CC^K)^\rho_\sigma\CC^\gamma_\rho =
[\delta_\rho(\CC^K)^\rho_\mu]\CC^\sigma_\nu\CC^\gamma_\rho -
C^\rho_{\nu\sigma}(\CC^K)^\sigma_\mu\CC^\gamma_\rho
$$
By the usual Leibniz rule for $\delta_\rho$, this yields
$$
C^\gamma_{\mu\nu}(\CC^K)^\rho_\sigma\CC^\gamma_\rho =
\delta_\rho (\CC^{K+1})^\gamma_\rho\CC^\sigma_\nu -
(\CC^K)^\rho_\mu C^\gamma_{\rho\sigma}\CC^\sigma_\nu -
C^\rho_{\nu\sigma}(\CC^K)^\sigma_\mu\CC^\gamma_\rho.
$$
The identity (\ref{eq:CCN}) follows for $N = K+1$ if the second and third
summand on the right hand side add up to 
$-C^\gamma_{\nu\sigma}(\CC^{K+1})^\sigma_\mu$.
After renaming the indices, one brings the sum of these two to the
form
$$
(\CC^K)^\rho_\mu (- C^\sigma_{\nu\lambda}C^\gamma_{\rho\sigma} +
C^\sigma_{\nu\rho}C^\gamma_{\lambda\sigma})\partial^\lambda =
-(C^K)^\rho_\mu C^\sigma_{\rho\lambda}\partial^\lambda
C^\gamma_{\nu\sigma} = -(C^{K+1})^\sigma_\mu C^\gamma_{\nu\sigma}
$$
as required. The Jacobi identity is used for the 
equality on the left.\end{pf}

\subsection*{A.2 Cofiltered vector spaces and completions}

We sketch the formalism treating the algebraic duals $U(\gg)^*$ and 
$S(\gg)^*$ of filtered algebras $U(\gg), S(\gg)$ as cofiltered algebras.
The reader can treat them alternatively as topological algebras:
the basis of neighborhoods of $0$ 
in the formal adic topology of $U(\gg)^*$ and $S(\gg)^*$ is
given by the annihilator ideals $\operatorname{Ann}\,U_i(\gg)$ and 
$\operatorname{Ann}\,S_i(\gg)$, consisting of functionals vanishing
on the $i$-th filtered component. A {\bf cofiltration} on a vector 
space $A$ is an inverse sequence of epimorphisms of its 
quotient spaces $\ldots \to A_{i+1}\to A_i\to A_{i-1}\to\ldots\to A_0$;
denoting the quotient maps
$\pi_i:A\to A_i$ and $\pi_{i,i+k}:A_{i+k}\to A_i$, 
the identities $\pi_i = \pi_{i,i+k}\circ\pi_{i+k}$,  
$\pi_{i,i+k+l} = \pi_{i,i+k}\circ\pi_{i+k,i+k+l}$ are required to hold. 
The limit $\underset\longleftarrow\lim{}_r A_r$
consists of {\it threads}, i.e. the
sequences $(a_r)_{r\in\mathbb{N}_0}\in\prod_r A_r$
of compatible elements, $a_r = \pi_{r,r+k}(a_{r+k})$.
The canonical map $A\to\hat{A}$ to the {\bf completion} 
$\hat{A} := \underset\longleftarrow\lim{}_r A_r$ is $1$-$1$ if
$\forall a\in A$ $\exists r\in\mathbb{N}_0$ such that $\pi_r(a)\neq 0$.
The cofiltration is {\bf complete} if the canonical map $A\to\hat{A}$
is an isomorphism. 
{\bf Strict morphisms} of cofiltered vector spaces $A\to B$ are 
the linear maps which induce the levelwise
maps $A_r\to B_r$ on the quotients.
(This makes the category of complete cofiltered vector spaces
more rigid than the category of pro-vector spaces.) 
We say that $a= (a_r)_r$ has the {\bf cofiltered degree} $\geq N$ 
if $a_r = 0$ for $r<N$.
In our main example, $U_i(\gg)^*:= (U(\gg)^*)_i := 
U(\gg)^*/\operatorname{Ann}\,U_i(\gg)\cong (U_i(\gg))^*$
and similarly for $S(\gg)^*\cong\hat{S}(\gg^*)$. 
We use lower indices both
for filtrations and for cofiltrations 
(but upper for gradations!). 
Given a family of elements in $A$,
$a:\Lambda\to A$, $\lambda\mapsto a_\lambda$,
the expression  ('abstract infinite sum')
$\sum_{\lambda\in\Lambda}a_\lambda$ is called a {\bf formal sum}
if for each $r\geq 0$, there is only finitely many $\lambda$ 
such that $\pi_r(a_\lambda)\neq 0$ hence 
$\pi_r(\sum_{\lambda\in\Lambda}a_\lambda)
:= \sum_{\lambda\in\Lambda}\pi_r(a_\lambda)\in A_r$ is well defined; and therefore
there is well defined thread $(\pi_r(\sum_{\lambda\in\Lambda}a_\lambda))_r\in\hat{A}$, the value of the formal sum.

The usual tensor product $A\otimes B$ of cofiltered vector spaces
is cofiltered with the $r$-th cofiltered component 
(see~\cite{StojicPhD})
\begin{equation}\label{eq:cofiltensor}
(A\otimes B)_r = 
\frac{A\otimes B}{\cap_{p+q=r}\operatorname{ker}\,\pi^A_p
\otimes\operatorname{ker}\,\pi^B_q}.
\end{equation}
$(A\otimes B)_r$ is an abelian group of 
finite sums of the form $\sum_\lambda a_\lambda\otimes b_\lambda\in A\otimes B$ 
modulo the additive relation of equivalence $\sim_r$ for which
$\sum a_\mu\otimes b_\mu \sim_r 0$ iff $\pi_p(a_\mu)\otimes\pi_q(b_\mu) = 0$
in $A_p\otimes B_q$ for all $p, q$ such that $p+q=r$. 
Define the {\bf completed tensor product} 
$A\hat\otimes B  = \underset\longleftarrow\lim{}_r (A\otimes B)_r$,
equipped with the same cofiltration, $(A\hat\otimes B)_r := (A\otimes B)_r$.
An element in $A\hat\otimes B$ is thus the class of equivalence of a
formal sum $\sum_\lambda a_\lambda\otimes b_\lambda$  
such that for any $p$ and $q$ there are at most
finitely many $\lambda$ such that 
$\pi^A_p(a_\lambda)\otimes\pi^B_q(b_\lambda)\neq 0$. 
Alternatively, we can equip  $A\otimes B$ with a bicofiltration 
($\mathbb{N}_0\times\mathbb{N}_0$-cofiltration), $(A\otimes B)_{r,s} =
A_r\otimes B_s$. Observe the inclusions
$\operatorname{ker}\,\pi_{r+s}\otimes \operatorname{ker}\,\pi_{r+s}
\subset \cap_{p+q=r+s}\operatorname{ker}\,\pi_p\otimes \operatorname{ker}\,\pi_q\subset\operatorname{ker}\,\pi_r\otimes \operatorname{ker}\,\pi_s$,
which induce projections 
$A_{r+s}\otimes B_{r+s}\twoheadrightarrow (A\otimes B)_{r+s}\twoheadrightarrow
A_r\otimes B_s$ for all $r,s$; by passing to the limit we see that 
the completion with respect to the bicofiltration
and with respect to the original cofiltration are equivalent (and alike
statement for the convergence of infinite sums inside $A\hat\otimes B$).
A linear map among cofiltered vector spaces is
{\bf distributive over formal sums} if it sends 
formal sums to formal sums summand by summand (formal
version of $\sigma$-additivity). This property is weaker than
being a strict morphism of complete cofiltered vector spaces.
In fact~(\cite{StojicPhD}), a linear map $f:C\to D$ 
is distributive over formal sums iff 
$\forall s$ $\exists r$ and a linear map $f_{sr}:C_r\to D_s$
such that $\pi_s\circ f = f_{sr}\circ\pi_r$ (in the strict case we required
$s=r$). If $A$ and $B$ are complete, we can also consider maps
$A\otimes B\to C$ distributive over formal sums in each argument
separately. Unlike the strict morphisms of cofiltered spaces, 
such a map does not need to extend to a map $A\hat\otimes B\to\hat{C}$
distributive over formal sums in $A\hat\otimes B$ 
(continuity in each argument separately does not imply 
the joint continuity).

A (strict) {\bf cofiltered algebra} $A$ 
(e.g.~$\hat{S}(\gg^*)$) is a monoid internal to 
the $\genfd$-linear category of complete cofiltered vector spaces, 
strict morphisms and with the tensor product 
$\hat\otimes$~(\cite{StojicPhD}). The bilinear associative unital
multiplication map $\hat{m}:A\hat\otimes A\to A$ is a strict morphism,
hence inducing linear maps $m_r:(A\otimes A)_r\to A_r$ for all $r$.
In other words, $A\hat\otimes A\ni\sum_\lambda a_\lambda\otimes b_\lambda
\overset{\hat{m}}\mapsto \sum_\lambda a_\lambda\cdot b_\lambda\in A$, where
$(\sum_\lambda a_\lambda\cdot b_\lambda)_r$ is an equivalence class in $A_r$
of $(\pi_r\circ\hat{m})(\sum'_\lambda a_\lambda\otimes b_\lambda)$, 
where $\sum'$ 
denotes the {\it finite sum} over all $\lambda$ 
such that $\exists (p,q)$ with $p+q=r$ and 
$\pi_p(a_\lambda)\otimes\pi_q(b_\lambda) \neq 0$. 

Any vector subspace $W$ of a cofiltered vector space $V$ is cofiltered
by $W_p := V_p\cap W$ with a canonical linear map 
$\underset\longleftarrow\lim\,W_p\to\underset\longleftarrow\lim\,V_p=\hat{V}$, 
whose image is a cofiltered subspace $\hat{W}_{\hat{V}}\subset\hat{V}$, 
the {\it completion} of $W$ in $\hat{V}$.
This is compatible with many additional structures, so defining 
the completions of sub(bi)modules and ideals 
(thus $\hat{I}$, $\hat{I}'$, $\hat{I}^{(r)}$, $\hat{I}'^{(r)}$,
$\hat{\bar{I}}$, $\hat{\bar{I}}'$, $\hat{\bar{I}}^{(r)}$, 
$\hat{\bar{I}}'^{(r)}$ 
in Sections~\ref{sec:bimodule} and~\ref{sec:bialgebroid}). 
If $U$ is an associative algebra, 
$A_U$ a right $U$-module and ${}_U B$ a left $U$-module, 
where both modules are cofiltered, then define $A\hat\otimes_U B$ 
as the quotient of $A\hat\otimes B$
by the completion of $\operatorname{ker}\,(A\otimes B\to A\otimes_U B)$ 
in $A\hat\otimes B$.

In this article, the completed tensor product 
$U(\gg^L)\hat\otimes\hat{S}(\gg^*)$ is defined by equipping
the filtered ring $U(\gg^L)$ with the {\it trivial cofiltration} 
$U(\gg^L)$, in which every cofiltered component is the entire $U(\gg)$
(and carries the discrete topology). 
The elements of $U(\gg^L)\hat\otimes\hat{S}(\gg^*)$ 
are given by the formal sums  $\sum u_\lambda\otimes a_\lambda$ 
such that $\forall r$, 
$\pi_r(a_\lambda) = 0$ for all but finitely many $\lambda$.
The basis of neighborhoods of $0$ in 
$U(\gg^L)\hat\otimes\hat{S}(\gg^*)$ 
consists of the subspaces 
$\genfd\hat{f}\otimes\prod_{p>r}S^p(\gg^*)$ 
for all $\hat{f}\in U(\gg)$ and $r\in\mathbb{N}$. 
The right Hopf action 
$a\otimes\hat{u}\mapsto\bm\phi_+(\hat{u})(a)$ {\it admits} a 
completed smash product:

\begin{theorem}\label{thm:completedsmash} 
The multiplication in $H^L = U(\gg^L)\sharp_{\bm\phi_+}\hat{S}(\gg^*)$
extends to a unique multiplication $\hat{m}$ on 
$U(\gg^L)\hat\otimes\hat{S}(\gg^*)$ which distributes 
over formal sums in each argument, 
forming the {\bf completed smash product algebra}
$\hat{H}^L = U(\gg^L)\hat\sharp\hat{S}(\gg^*)$. Likewise, the
multiplication on  $H^R = \hat{S}(\gg^*)\sharp_{\bm\tilde\phi_-}U(\gg^R)$
extends to $\hat{S}(\gg^*)\hat\otimes U(\gg^R)$ forming
$\hat{H}^R = \hat{S}(\gg^*)\hat\sharp U(\gg^R)$.
However, there are no cofiltered algebra structures on~$\hat{H}^L$, 
because the multiplication does not
distribute over formal sums in $H^L\hat\otimes H^L$.
\end{theorem} 
\begin{pf}  The extended multiplication is well defined by a formal sum 
$\sum_{\lambda,\mu} (u_\lambda\sharp a_\lambda)(u'_\mu\sharp a'_\mu)
= \sum_{\lambda,\mu} u_\lambda u'_{\mu(1)}\sharp 
 \bm\phi_+(u'_{\mu(2)})(a_\lambda) a'_\mu$ 
if for all $r\in\mathbb{N}_0$
the number of pairs $(\mu,\lambda)$ such that 
$u_\lambda u'_{\mu(1)}\otimes\pi_r(\bm\phi_+(u'_{\mu(2)})(a_\lambda)
\cdot a'_\mu)\neq 0$
(only Sweedler summation) is finite. There are only
finitely many $\mu$ such that $\pi_r(a'_\mu)\neq 0$; 
only those contribute to the sum because $\pi_k(a)\pi_l(b)=0$ 
implies $\pi_{k+l}(ab)=0$ in any cofiltered ring. 
For each such $\mu$ fix a representation of 
$\Delta(u_\mu)$ as a finite sum 
$\sum_k u_{\mu(1)k}\otimes u_{\mu(2)k}$ and denote by $K(\mu)$
the maximal over $k$ filtered degree of $u_{\mu(2)k}$ and by $L(\mu)$
the minimal cofiltered degree of $a'_\mu$. 
By Lemma~\ref{lemphi} (iii) and induction we see that if 
$a_\lambda\in\hat{S}(\gg^*)_s$ then 
$\bm\phi_+(u_{\mu(2)k})(a_\lambda)\in\hat{S}(\gg^*)_{s-K(\mu)}$. Hence for
each $\mu$ 
there are only finitely many $\lambda$ for which $s-K(\mu)+L(\mu)\leq r$.
That is sufficient for the conclusion.
More details will be exhibited elsewhere.   
\end{pf}

\footnotesize{

}
\end{document}